\documentclass{article}

\usepackage[microtype]{preamble}
\usepackage[margin=1in]{geometry}
\usepackage{authblk}

\usepgfplotslibrary{patchplots}

\title{Adaptive sampling-based quadrature rules for efficient Bayesian~prediction}
\author[1,2]{L.\,M.\,M.~van~den~Bos}
\author[1]{B.~Sanderse}
\author[2]{W.\,A.\,A.\,M.~Bierbooms}

\affil[1]{Centrum~Wiskunde~\&~Informatica, P.O.~Box 94079, 1090~GB, Amsterdam}
\affil[2]{Delft~University~of~Technology, P.O.~Box 5, 2600~AA, Delft}

\newcommand{\x}{\ensuremath\mathbf{x}}
\newcommand{\y}{\ensuremath\mathbf{y}}
\newcommand{\z}{\ensuremath\mathbf{z}}
\newcommand{\A}{\ensuremath\mathcal{A}}
\newcommand{\I}{\ensuremath\mathcal{I}}

\DeclareMathOperator*{\Cov}{Cov}

\DeclareMathOperator*{\linspan}{span}

\begin{document}
\maketitle

\begin{abstract}
\noindent A novel method is proposed to infer Bayesian predictions of computationally expensive models. The method is based on the construction of quadrature rules, which are well-suited for approximating the weighted integrals occurring in Bayesian prediction. The novel idea is to construct a sequence of nested quadrature rules with positive weights that converge to a quadrature rule that is weighted with respect to the posterior. The quadrature rules are constructed using a proposal distribution that is determined by means of nearest neighbor interpolation of all available evaluations of the posterior. It is demonstrated both theoretically and numerically that this approach yields accurate estimates of the integrals involved in Bayesian prediction. The applicability of the approach for a fluid dynamics test case is demonstrated by inferring accurate predictions of the transonic flow over the RAE2822 airfoil with a small number of model evaluations. Here, the closure coefficients of the Spalart--Allmaras turbulence model are considered to be uncertain and are calibrated using wind tunnel measurements.
\end{abstract}

\begingroup
\endgroup

\section{Introduction}
Computer simulation with uncertainties often requires calculating the expectation with respect to a probability distribution of a complex and computationally costly function. Various approaches exist to numerically estimate such integrals, for example Monte Carlo approaches~\cite{Caflisch1998} or deterministic quadrature rule approaches~\cite{Brass2011}. In these approaches it is often assumed that samples from the probability distribution can be constructed straightforwardly or are readily available. However, this is not always the case, for instance in the case of Bayesian calibration problems~\cite{Kennedy2001,Gelman2013} the probability distribution encompasses the computationally costly model and is highly non-trivial. This is known as Bayesian prediction and is the main topic of this article.

A commonly used approach to calculate weighted integrals with respect to such non-trivial distributions is importance sampling~\cite{Tokdar2009,Derflinger2010,Botts2011}. In this approach samples are drawn from a \emph{proposal distribution} and by means of weighting these samples are used to determine integrals with respect to a \emph{target distribution}, i.e.\ the distribution for which it is difficult to construct samples. If the proposal distribution is chosen well, the constructed samples yield similar convergence rates as samples constructed directly with the target distribution. A major disadvantage is that constructing the proposal distribution is not trivial and often requires a priori knowledge about the target distribution. A well-known extension of importance sampling to a Bayesian framework is formed by Markov chain Monte Carlo methods~\cite{Metropolis1953,Hastings1970}, where the proposal distribution is chosen such that the samples form a Markov chain. In this case the proposal distribution is not constructed beforehand, though a large number of costly evaluations of the posterior is typically necessary to construct the sequence of samples.

There exist various approaches to improve the performance of Monte Carlo techniques in Bayesian inference. For example, using preconditioning based on Laplace's method~\mbox{\cite{Wacker2017}}, where an approximation of the posterior is constructed using the maximum a-posteriori estimate and a Gaussian distribution; using transport maps~\mbox{\cite{Peherstorfer2019}}, where a low-fidelity model is used to precondition a Markov chain Monte Carlo sampler; or by adaptively tuning the proposal distribution~\mbox{\cite{Vrugt2009}}. However, many preconditioning methods are significantly less efficient if the posterior cannot be represented well using a Gaussian distribution. Techniques from the field of machine learning can also be applied~\mbox{\cite{Ahn2012}}, though the focus of those techniques is mostly on the treatment of large data sets.

Even if samples can be drawn from the posterior (possibly by using any of these improvements), the error of estimating an integral by means of random sampling decays as $\mathcal{O}(1/\sqrt{K})$ (with $K$ the number of samples). For example, to double the accuracy of the estimate the number of samples must be quadrupled. Therefore using importance sampling techniques, or Monte Carlo techniques in general, in conjunction with a computationally costly model is often infeasible. In summary, estimating integrals by directly sampling from the posterior is often difficult or very expensive with existing methods.

A commonly used approach to alleviate the high computational cost of sampling and the slow convergence is to construct a surrogate (or ``response surface'') and use the surrogate instead of the complex model for inferring the relevant statistics. Well-known approaches include Gaussian process regression~\cite{Stuart2017,Bilionis2013} and polynomial approximation~\cite{Marzouk2009,Birolleau2014,Schwab2012,Schillings2013,Bos2018}. If the weighting distribution is not known explicitly, it is difficult to construct an optimal surrogate (and moreover it often requires analyticity or smoothness), so often the surrogate is constructed such that it is globally accurate (or the distribution is tempered~\cite{Li2014,Bos2018}, which yields a similarly accurate surrogate). It has been demonstrated rigorously that a posterior determined with a globally accurate surrogate converges to the true posterior~\cite{Yan2017,Butler2018a}.

At the same time, in the context of Bayesian prediction the interest is not directly in accurately approximating the \emph{model}, but in calculating \emph{moments} of the model with respect to the posterior. These moments are integral quantities and there exist many efficient numerical integration techniques for this purpose that do not explicitly construct a surrogate. In this article, the focus is specifically on quadrature rules due to their high accuracy and flexibility in the location of the nodes. Existing quadrature rule approaches~\cite{Brass2011} provide rapid converging estimates of the integral for smooth functions, without suffering from a deteriorated convergence rate if the integrand is not smooth (as often observed in interpolation approaches). However, the state-of-the-art quadrature rule approaches such as Gaussian quadrature rules~\cite{Golub1969}, multivariate sparse grids~\cite{Novak1999,Smolyak1963,Hewitt2019}, or heuristic optimization approaches~\cite{Jakeman2017,Narayan2014} also require much knowledge about the distribution under consideration, so these are not directly applicable to the problems of interest discussed here.

Therefore, the main goal of this work is to propose a novel flexible method for Bayesian prediction that combines the high convergence rates of quadrature rules with the flexibility of adaptive importance sampling. The idea is to construct a sequence of quadrature rules and a sequence of proposal distributions, where the proposal distributions are constructed using all available model evaluations. The sequence of quadrature rules converges to a quadrature rule that incorporates the (expensive) probability distribution (i.e.\ the posterior). By ensuring that all constructed quadrature rules have high degree and have positive weights, high convergence rates are obtained (depending on the smoothness of the integrand). The rules are specifically tailored to determine moments weighted with a posterior and require a small number of nodes to obtain an accurate estimation. Moreover it can be demonstrated rigorously that estimations made with the sequence of quadrature rules converge to the exact integral for a large class of integrands and distributions.

This article is structured as follows. In Section~\ref{sec:problemsetting} the problem of Bayesian prediction with a quadrature rule is briefly introduced, with a main focus on quadrature rule methods and their relevant properties. In Section~\ref{sec:quadraturerule} our main contribution is discussed: a method to construct a quadrature rule that converges to a quadrature rule of the posterior. The method is analyzed theoretically in Section~\ref{sec:erranalysis} and to demonstrate the applicability of the quadrature rule to computational problems, the technique is applied to some test functions in Section~\ref{sec:numerics}. Moreover the calibration of the turbulence closure coefficients of the RAE2822 airfoil is discussed. The article is summarized and concluded in Section~\ref{sec:conclusion}.

\section{Preliminaries}
\label{sec:problemsetting}

\subsection{Bayesian prediction}
The key goal of this article is to propose an efficient method to infer Bayesian predictions involving a computationally expensive model. However, the setting where the proposed methods can be applied is slightly more general than that. Let $u: \Omega \to \mathbb{R}^n$ be a continuous function with $\Omega \subset \mathbb{R}^d$ ($d = 1, 2, 3, \dots$) that describes the quantity of interest and let $\rho: \Omega \to [0, \infty)$ be a continuous and bounded probability density function (PDF). Throughout this work, it is assumed that $u$ can be evaluated exactly without any significant numerical error. Even though this is not a realistic assumption in practice, the stability of the methods proposed in this work guarantees that any numerical error present in the computational model does not affect the accuracy of the proposed approaches (i.e.\ errors ``do not blow up''). Consequently determining Bayesian predictions can be described as determining the expectation of a function $f: \Omega \to \mathbb{R}$, with $f(\mathbf{X}) = F(u(\mathbf{X}))$ where $\mathbf{X}$ is a $d$-variate random vector with PDF $\rho$ (so $F$ is a function defined on the \emph{range} of $u$). This integral will be denoted by the operator $\I^\rho$ throughout this article, i.e.
\begin{equation}
	\label{eq:leintegral}
	\I^\rho f \coloneqq \mathbb{E}[f(\mathbf{X})] = \int_\Omega f(\x) \, \rho(\x) \dd \x = \int_\Omega F(u(\x)) \, \rho(\x) \dd \x.
\end{equation}

Integrals of this type have been studied many times in the framework of uncertainty propagation and many efficient methods exist to approximate them~\cite{Xiu2010,Sullivan2015,LeMaitre2010}. The approach taken in this work is based on a nearest neighbor interpolant of $\rho(\x)$, which is used to construct quadrature rules for approximating the integral of \eqref{eq:leintegral}. Geometric approaches based on nearest neighbor interpolation and the closely related Voronoi tessellation have been used succesfully to construct estimates of integrals as considered in this article~\mbox{\cite{Guessab2008,Luo2009,Butler2017,Du1999}}. However, in these cases the distribution $\rho$ is known explicitly or is computationally tractable beforehand (possibly up to constant). In this article, the computationally challenging part is that $\rho(\x)$ is obtained by Bayesian analysis, i.e.\ it is, up to a constant, a \emph{posterior}: 
\begin{equation}
	\rho(\x) \coloneqq q(\z \mid \x) \, q(\x),
\end{equation}
where $\z$, $q(\z \mid \x)$, and $q(\x)$ are the measurement data, the likelihood, and the prior respectively. The posterior, denoted by $q(\x \mid \z)$, is obtained by scaling $\rho(\x)$ such that it is a distribution:
\begin{equation}
	q(\x \mid \z) = \frac{q(\z \mid \x) \, q(\x)}{\int_\Omega q(\z \mid \x) \, q(\x) \dd \x} = \frac{\rho(\x)}{\int_\Omega \rho(\x) \dd \x}.
\end{equation}
The scaling factor (often called the \emph{evidence}) that is present in this definition of $\rho$ is neglected throughout this article, as it is not of importance for the proposed method. It will be shown that the method proposed in this work is stable (see Section~\ref{subsubsec:degreepos}) regardless of the exact value of the evidence. Moreover, this is also demonstrated numerically, by considering an example of a case with small evidence (see Section~\ref{subsubsec:peakypeaky}).

Bayesian prediction consists of assessing the \emph{posterior predictive distribution}, which describes the probability of observing new data $\widehat{\z}$ given the existing data $\z$. It is obtained by marginalizing over the posterior as follows~\cite{Edeling2014-2}:
\begin{equation}
	q(\widehat{\z} \mid \z) = \int_\Omega q(\widehat{\z}, \x \mid \z) \dd \x = \int_\Omega q(\widehat{\z} \mid \x, \z) \, q(\x \mid \z) \dd \x = \int_\Omega q(\widehat{\z} \mid \x) \, q(\x \mid \z) \dd \x.
\end{equation}
Here, we assume that $\widehat{\z}$ and $\z$ are conditionally independent given $\x$. Moments of the posterior predictive distribution can be used among others to obtain point predictions. These moments can often be expressed explicitly as integrals over the computational model $u$.

In general we require throughout this article that the distribution $\rho$ is such that an evaluation of $\rho(\x)$ requires an evaluation of $u(\x)$. Moreover if $u$ has been evaluated for a specific value of $\x$, the value of $\rho(\x)$ can be determined without any significant computational cost. This is often the case if $\rho(\x)$ forms a posterior in a Bayesian framework and this property is the key observation that will be leveraged to approximate \eqref{eq:leintegral}.

A well-known example of a prior and a likelihood for a scalar model $u$ are a uniform prior, defined as $q(\x) \propto 1$ for $\x \in \Omega$, and a Gaussian likelihood, defined as follows for given standard deviation $\sigma$:
\begin{equation}
	q(\z \mid \x) \propto \exp\left[ -\frac{1}{2} \frac{\sum_{k=1}^m (u(\x) - z_k)^2}{\sigma^2} \right], \text{ with measurement data $\z = \trans{(z_1, \dots, z_m)}$ known.}
\end{equation}
In this particular case, the expectation of the posterior predictive distribution can be explicitly expressed as follows:
\begin{multline}
	\label{eq:prediction}
	\mathbb{E}[\widehat{\z} \mid \z] = \int_Z \widehat{\z} \, q(\widehat{\z} \mid \z) \dd \widehat{\z} = \int_\Omega \int_Z \widehat{\z} \, q(\widehat{\z} \mid \x) \, q(\x \mid \z) \dd \widehat{\z} \dd \x \\
	= \int_\Omega q(\x \mid \z) \int_Z \widehat{\z} \, q(\widehat{\z} \mid \x) \dd \widehat{\z} \dd \x = \int_\Omega q(\x \mid \z) \, \mathbb{E}[\widehat{\z} \mid \x] \dd \x = \int_\Omega q(\x \mid \z) \, u(\x) \dd \x.
\end{multline}
Here, $Z$ denotes the space that contains all possible values of the data. It is not necessary to derive this space formally, as it is only used in an intermediate step of the derivation. Notice that the obtained integral is the expectation of the computational model with the parameters distributed according to the posterior. Various statistical models that combine measurement error and model error exist for the purpose of Bayesian model calibration. We refer the interested reader to~\citet{Kennedy2001} and will consider such an extensive model in a numerical example discussed in Section~\ref{sec:numerics}. Vector-valued $u$ can be incorporated straightforwardly in the framework of this article, since only the posterior distribution $q(\x \mid \z)$ is used to construct numerical integration techniques.

Throughout this article, two assumptions are imposed on the statistical setting. Firstly, we assume that the statistical moments of the prior are finite, i.e.
\begin{equation}
	\int_\Omega \varphi(\x) \, q(\x) \dd \x < \infty, \text{ for all polynomials $\varphi$}.
\end{equation}
This implies among others that the prior is not improper. Secondly, it is assumed that the posterior is continuous and bounded. These two assumptions are not restrictive, e.g.\ statistical models with Gaussian random variables fit in this setting. Moreover, these assumptions are only necessary for the theoretical analysis of the approach discussed in this article. The approach requires only straightforward modifications for applying it to cases where an improper prior or discontinuous posterior is considered, and an example of the latter is discussed in more detail in Section~\ref{subsec:explicit}.

\subsection{Quadrature rules}
\label{subsec:quadrules}
The integral from \eqref{eq:leintegral} is approximated by means of a quadrature rule, consisting of nodes $\{\x_k\}_{k=0}^N \subset \Omega$ and weights $\{w_k\}_{k=0}^N \subset \mathbb{R}$:
\begin{equation}
	\label{eq:leapproximation}
	\int_\Omega f(\x) \, \rho(\x) \dd \x \approx \sum_{k=0}^N w_k f(\x_k).
\end{equation}
Notice that the distribution $\rho(\x)$ does not appear explicitly on the right hand side of this approximation, but that it is implicitly incorporated in the values of the nodes $\x_k$ and the weights $w_k$. We will denote the quadrature rule approximation by means of the linear operator $\A^\rho_N$, i.e.
\begin{equation}
	\A^\rho_N f \coloneqq \sum_{k=0}^N w_k f(\x_k).
\end{equation}
If $\rho(\x) = q(\z \mid \x) \, q(\x)$, a quadrature rule with respect to the posterior $q(\x \mid \z)$ can be obtained by rescaling the weights of a quadrature rule of $\rho$ such that $\sum_{k=0}^N w_k = 1$, provided that $\A^\rho_N 1 = \I^\rho 1$.

Three properties are relevant for quadrature rules that are applied to computationally expensive models: the exactness of the quadrature rule, positivity of the weights, and nesting. The first two properties ensure numerical stability of the quadrature rule and make that the quadrature rule converges. These two properties and their relevance are discussed in Section~\ref{subsubsec:degreepos}. Nesting allows for straightforward refinement of the quadrature rule approximation, and is briefly discussed in Section~\ref{subsubsec:nesting}.

\subsubsection{Exactness and positive weights}
\label{subsubsec:degreepos}
Often the exactness of a quadrature rule is related to the degree of the polynomials it integrates exactly. We choose to define the exactness slightly more general to facilitate refinements of arbitrary numbers of nodes.

To this end, let $\varphi_0, \dots, \varphi_D$ be $D+1$ $d$-variate polynomials defined on $\Omega$ and consider the space $\Phi_D = \linspan\{ \varphi_0, \dots, \varphi_D \}$. Then the goal is to construct a quadrature rule such that it integrates all functions $\varphi \in \Phi_D$ exactly. This is, due to linearity, equivalent to
\begin{equation}
	\A^\rho_N \varphi_j = \I^\rho \varphi_j \text{ for } j = 0, \dots, D.
\end{equation}
Without loss of generality, we assume throughout this article that the polynomials are sorted graded lexicographically. This ensures that referring to $\varphi_k$ for any $k$ is well-defined and that $\Phi_D \subset \Phi_{D+1}$ for any $D$, but other monomial orders that are a well-order can be used straightforwardly. Moreover we assume that $\varphi_0$ is the constant polynomial, as that allows to straightforwardly scale the quadrature rule weights to incorporate the evidence. In the univariate case, this implies that $D$ is the degree of the quadrature rule, since $\Phi_D$ contains all polynomials of degree $D$ and less. In the multivariate case, the degree of the quadrature rule equals $Q$ if $\dim \Phi_D = \dim \mathbb{P}(Q, d)$ where $\mathbb{P}(Q, d)$ denotes all $d$-variate polynomials of a degree $Q$.

The close relation between the nodes and the weights of a quadrature rule is described by the $(D+1) \times (N+1)$ \emph{Vandermonde matrix} $V_D(X_N)$:
\begin{equation}
	\underbrace{\begin{pmatrix}
		\varphi_0(\x_0) & \cdots & \varphi_0(\x_N) \\
		\vdots & \ddots & \vdots \\
		\varphi_D(\x_0) & \cdots & \varphi_D(\x_N)
	\end{pmatrix}}_{V_D(X_N)}
	\begin{pmatrix}
		w_0 \\
		\vdots \\
		w_N
	\end{pmatrix}
	=
	\begin{pmatrix}
		\mu_0 \\
		\vdots \\
		\mu_D
	\end{pmatrix},
\end{equation}
where $\mu_j = \I^\rho \varphi_j$. If $\varphi_j$ are monomials, $\mu_j$ are the raw moments of the distribution $\rho$. If $N=D$ the quadrature rule is called \emph{interpolatory} and the linear system can be used to determine the weights given the nodes, provided that the Vandermonde matrix is non-singular. In the univariate case, this holds if and only if the nodes are distinct. In the multivariate case, it is less straightforward to ensure that the Vandermonde matrix is non-singular, but all quadrature rules constructed in this work have a non-singular Vandermonde matrix.

If the weights are non-negative (i.e.\ $w_k \geq 0$ for all $k$), the dimension of $\Phi_D$ is a measure for the accuracy of the quadrature rule and increasing $D$ yields a converging estimate for sufficiently smooth functions. To see this, assume that $\Omega$ is bounded and let $\widehat{\varphi}$ be the best approximation of the function $u$ in $\Phi_D$ measured in the $\infty$-norm (assuming without loss of generality that it exists):
\begin{equation}
	\widehat{\varphi} \coloneqq \argmin_{\varphi \in \Phi_D} \| u - \varphi \|_\infty, \text{ with } \| u - \varphi \|_\infty \coloneqq \sup_{\x \in \Omega} |u(\x) - \varphi(\x)|.
\end{equation}
Then:
\begin{equation}
	\label{eq:lebesgue}
	| \I^\rho u - \A^\rho_N u | \leq \| \I^\rho \|_\infty \| u - \widehat{\varphi} \|_\infty + \| \A^\rho_N \|_\infty \| u - \widehat{\varphi} \|_\infty = \left( \| \I^\rho \|_\infty + \| \A^\rho_N \|_\infty \right) \inf_{\varphi \in \Phi_D} \| u - \varphi \|_\infty,
\end{equation}
where
\begin{equation}
	\| \I^\rho \|_\infty = \int_\Omega 1 \rho(\x) \dd \x = \sum_{k=0}^N w_k \leq \sum_{k=0}^N |w_k| = \| \A^\rho_N \|_\infty.
\end{equation}
The inequality \eqref{eq:lebesgue} is commonly known as the Lebesgue inequality, see e.g.~\citet[Theorem~3.1.1]{Brass2011}. It is arguably best known for its usage in polynomial interpolation~\mbox{\cite{Ibrahimoglu2016,Brutman1996}}, but it straightforwardly carries over to the setting of quadrature rules (then the operator norm $\| \A_N^\rho \|_\infty$ is the Lebesgue constant).

If all weights are positive, then $|w_k| = w_k$ and therefore $\| \A^\rho_N \|_\infty = \| \I^\rho \|_\infty$. Hence
\begin{equation}
	| \I^\rho u - \A^\rho_N u | \leq 2 \| \I^\rho \|_\infty \inf_{\varphi \in \Phi_D} \| u - \varphi \|_\infty.
\end{equation}
As explained before, the weights of the quadrature rule can be scaled such that $\rho$ forms a probability density function. In that case, it holds that $\| \I^\rho \|_\infty = \| \A^\rho_N \|_\infty = 1$ and therefore \eqref{eq:lebesgue} simplifies to
\begin{equation}
	| \I^\rho u - \A^\rho_N u | \leq 2 \inf_{\varphi \in \Phi_D} \| u - \varphi \|_\infty.
\end{equation}

It is well-known that the right hand side of this inequality decays for increasing $D$ if $u$ can be approximated by a polynomial. This is among others the case if $u$ is absolutely continuous and in general the rate of decay of the right hand side of this inequality depends on the smoothness of $u$. If the function $u$ is univariate and $u \in C^r$, i.e.\ $u$ is $r$ times differentiable, the error usually decays with rate $r$, or in other words:
\begin{equation}
	\|u - \varphi \|_\infty \leq A D^{-r},
\end{equation}
where $A$ is a constant that solely depends on $u$. If $u$ is smooth (i.e.\ infinitely many times continuously differentiable), the error decays exponentially. The exact rate of decay can be derived using Jackon's inequality~\mbox{\cite{Jackson1982}}, but many more results on this topic exist. We refer the interested reader to~\citet{Watson1980} and the references therein for more information. The upper bound of the Lebesgue inequality is not sharp, as for unbounded $\Omega$ and discontinuous $u$ the estimate of a quadrature rule might converge, even though $\inf_\varphi \| u - \varphi \|_\infty$ is unbounded. It is out of the scope of this article to fully discuss all convergence properties of quadrature rules (but see e.g.~\citet{Brandolini2014} or~\citet{Brass2011}).

Notice that the bound described by the Lebesgue inequality depends on the dimension of $\Omega$. For multivariate $u$, the quantity $\| u - \varphi\|_\infty$ decays significantly slower. In general, it holds that if $u$ is $r$ times differentiable, the error decays at least with rate $r / d$, where $d$ is the dimension of the space $\Omega$. Estimating the decay of this quantity in multivariate space is significantly less straightforward than estimating it in univariate spaces, since computing the best approximation becomes significantly more difficult. Nonetheless, the convergence rates as stated here are sufficient for the analysis conducted in this work.

A quadrature rule with positive weights and $\sum w_k = 1$ is stable regardless of the computational and statistical model under consideration~\cite{Brass2011}. In other words, small variations in the integrand do not affect the accuracy of the rule as a whole. In particular, if $u$ is perturbed by a small error $\vareps > 0$, say $\tilde{u} = u \pm \vareps$, then
\begin{equation}
	| \A^\rho_N u - \A^\rho_N \tilde{u} | \leq \| \A^\rho_N \|_\infty \| u - \tilde{u} \|_\infty \leq \vareps.
\end{equation}
This stability result is particularly important for the usage of quadrature rules in a Bayesian setting, as it ensures that the rule yields stable estimates regardless of the exact value of the evidence (recall that the evidence is the scaling factor neglected in the definition of $\rho$). However, a small evidence will result in slow convergence, as will be demonstrated numerically in Section~\ref{subsubsec:peakypeaky}.

\subsubsection{Nesting}
\label{subsubsec:nesting}
Nesting means that smaller sets of quadrature rule nodes are contained in larger sets of quadrature rule nodes, i.e.\ $X_{N_1} \subset X_{N_2}$ for $N_1 < N_2$. Strictly speaking this forms a \emph{sequence} of quadrature rules, but we will call this with a little abuse of nomenclature simply a nested quadrature rule.

A nested quadrature rule allows for straightforward refinement of the quadrature rule approximation, as computationally costly function evaluations are being reused upon considering a larger quadrature rule. Moreover, if the weights of the quadrature rule are non-negative, it provides a straightforward way to assess the accuracy of a quadrature rule approximation by comparing two consecutive quadrature rule levels. This can be seen as follows:
\begin{equation}
	| \A^\rho_{N_1} u - \A^\rho_{N_2} u | \leq | \A^\rho_{N_1} u - \I^\rho u | + | \A^\rho_{N_2} u - \I^\rho u | \leq 4 \| \I^\rho \|_\infty \inf_{\varphi \in \Phi_D} \| u - \varphi \|_\infty,
\end{equation}
where $D$ is such that $\A^\rho_{N_1} \varphi = \A^\rho_{N_2} \varphi = \I^\rho \varphi$ for all $\varphi \in \Phi_D$. The difference between two consecutive quadrature rules also forms an upper bound:
\begin{equation}
	| \A^\rho_{N_1} u - \A^\rho_{N_2} u | \geq | \A^\rho_{N_1} u - \I^\rho u | - | \A^\rho_{N_2} u - \I^\rho u |,
\end{equation}
which is the difference between the errors of the two quadrature rules. This difference vanishes for $N_1, N_2 \to \infty$ if a converging sequence of quadrature rule is used.

\section{Bayesian predictions with an adaptive quadrature rule}
\label{sec:quadraturerule}
\begin{figure}
	\centering
	\includegraphics[width=\textwidth]{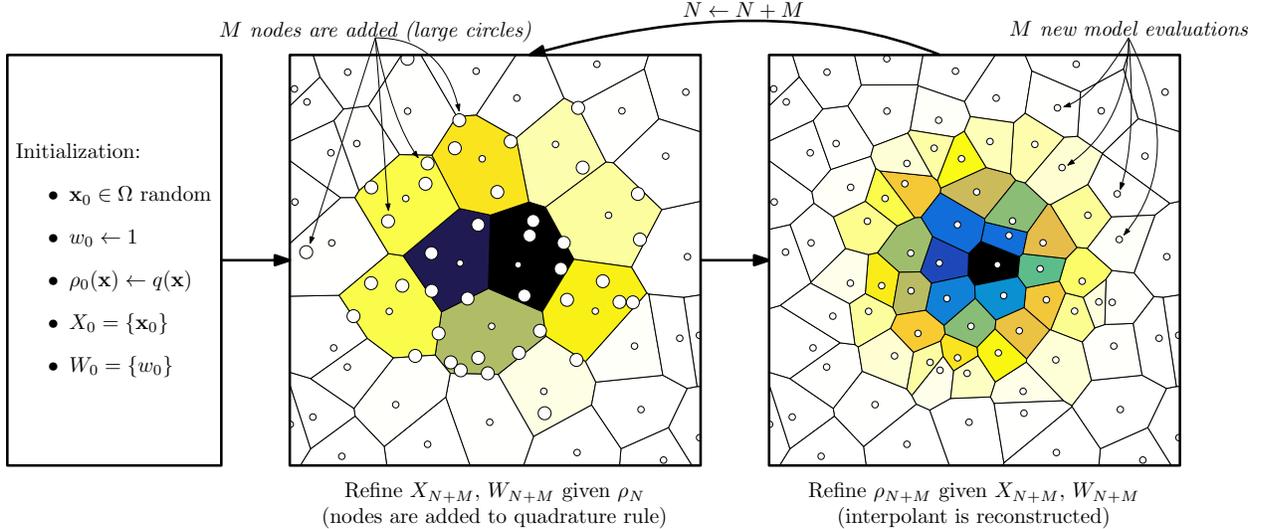}
	\caption{Schematic overview of the two main steps of the method to construct the quadrature rules proposed in this article.}
	\label{fig:algorithm}
\end{figure}
The main problem addressed in this article is accurate numerical estimation of Bayesian predictions. To this end, the main interest is the construction of quadrature rules to approximate integrals weighted with the posterior, a problem that was summarized in \eqref{eq:leapproximation} as follows:
\begin{equation}
	\I^\rho f = \int_\Omega f(\x) \, \rho(\x) \dd \x \approx \sum_{k=0}^N w_k f(\x_k) = \A^\rho_N f.
\end{equation}
Here, both $f$ and $\rho$ are functions of $u$, denoted by $f(\x) = F(u(\x))$ and $\rho(\x) = G(u(\x))$ for suitable $F$ and $G$. The functions $F$ and $G$ are such that they can be evaluated efficiently and quickly.

The focus of this article is mainly on using a quadrature rule to accurately infer Bayesian predictions, but various other statistical quantities that are not considered in this article can be inferred straightforwardly without additional model evaluations. For example, moments of the posterior distribution can be determined by integrating suitably chosen polynomials or an empirical cumulative distribution function can be constructed by considering the discrete probability density function $p(\x_k) = w_k$ for all $k$ (which is well-defined because all weights are non-negative).

\subsection{Constructing an adaptive quadrature rule}
The proposed procedure consists of iteratively refining a pair consisting of a quadrature rule and a distribution. A typical iteration consists of two steps. Firstly, a quadrature rule based on the current estimate of the distribution $\rho$ is constructed. Secondly, the estimate of the distribution is refined using all available model evaluations of $\rho$. The quadrature rules are constructed such that a sequence of nested rules is obtained. These steps are illustrated in Figure~\ref{fig:algorithm}.

More specifically, at the start of the $i$-th iteration a distribution $\rho_N$ is given, accompanied by evaluations of the model $u$ and distribution $\rho$ at nodes $\x_0, \dots, \x_N$. The first step consists of constructing a new quadrature rule with $N+M$ nodes that incorporates $\rho_N$ and reuses the available model evaluations, i.e.\ it is nested. Hence the goal is to determine $\x_{N+1}, \dots, \x_{N+M}$ and $w_0, \dots, w_{N+M}$ with $M$ as small as possible such that
\begin{equation}
	\sum_{k=0}^{N+M} w_k \, \varphi_j(\x_k) = \int_\Omega \varphi_j(\x) \, \rho_N(\x) \dd \x, \text{ for all $j = 0, \dots, D_i$}.
\end{equation}
Here, $D_i$ is a free parameter which defines the exactness of the rule, i.e.\ the dimension of $\Phi_{D_i}$. The second step of the iteration consists of constructing a refined approximation of the distribution $\rho$ using all available evaluations of the model $u$ and distribution $\rho$ (including the $M$ nodes that have been added this iteration). In this work, this is done by nearest neighbor interpolation. Consequently $\rho_{N+M}$ is obtained, which is used for the next iteration. The procedure can be initialized using a quadrature rule consisting of a single random node and a proposal distribution $\rho_0(\x) \equiv q(\x)$, i.e.\ by using the prior. The convergence can each iteration be assessed by comparing the current quadrature rule estimate with that of a previous iteration. The algorithm is sketched in Figure~\ref{fig:algorithm}.

The only free parameter in this construction is $D_i$, which can be chosen heuristically. For convergence, it is essential to enforce that $D_i \to \infty$ for $i \to \infty$. In this work, we will mainly use two variants: linear growth, i.e.\ $D_i = \mathcal{O}(i)$, and exponential growth, i.e.\ $D_i = \mathcal{O}(2^i)$.

The distribution $\rho_N$ is constructed by interpolation of all available point evaluations of $\rho$. It can be interpreted as a \emph{proposal distribution}, as often used in importance sampling strategies. The efficiency of the approach can be demonstrated mathematically by reusing techniques from importance sampling.

The construction of a quadrature rule that incorporates the distribution $\rho_N$ is non-trivial. We will employ the recently introduced implicit quadrature rule~\cite{Bos2018b} for this purpose, though we emphasize that any methodology to construct quadrature rules for a given distribution can be used in our approach. The implicit quadrature rule is a nested quadrature rule with positive weights that is constructed solely using samples from the distribution $\rho_N$. We discuss its construction in Section~\ref{subsec:implquadrule}. Given the quadrature rule nodes, the distribution that is used to refine the quadrature rule is constructed using nearest neighbor interpolation, which is explained in Section~\ref{subsec:interpolation}.

\subsection{Implicit quadrature rule}
\label{subsec:implquadrule}
The implicit quadrature rule~\cite{Bos2018b} is a quadrature rule constructed using samples from a distribution. Hence let $K+1$ samples $\y_0, \dots, \y_K$ be given, denoted as the set $Y_K$. In the iterative procedure considered in this work (i.e.\ Figure~\ref{fig:algorithm}), these are a large number of samples from $\rho_N$. Moreover, some model evaluations have been performed in previous iterations, so let $X_N \subset \Omega$ be $N+1$ nodes in $\Omega$ where the value of $u$ (and therefore $\rho$) is known. These nodes are denoted by $\x_0, \dots, \x_N$. Notice that $\rho_N$ is typically obtained by nearest neighbor interpolation of these nodes, which we will discuss in more detail in Section~\ref{subsec:interpolation}.

The goal is to find for given $D$ (we omit the subscript $i$ in this section) a subset of $M$ samples from $Y_K$, say $\x_{N+1}, \dots, \x_{N+M}$, with $M \ll K$ such that
\begin{equation}
	\label{eq:keyproperty}
	\sum_{k=0}^{N+M} w_k \, \varphi_j(\x_k) = \frac{1}{K+1} \sum_{k=0}^K \varphi_j(\y_k), \text{ for $j = 0, \dots, D$}.
\end{equation}
The rule is constructed such that the weights $w_k$ are non-negative, i.e.\ $w_k \geq 0$ (for all $k$). This ensures that the quadrature rule error is bounded, which can be seen by applying the Lebesgue inequality from \eqref{eq:lebesgue}.

Notice that only samples from $\rho_N$ are necessary to construct this quadrature rule, from which the approximate moments are computed. This is an advantage, since in higher dimensional spaces it is often computationally costly to exactly determine the moments of $\rho_N$.

The nodes $\x_0, \dots, \x_N$ do generally not form a subset of the samples $Y_K$. Ideally we would like to have that $M = 0$, which would imply that no new costly model evaluations are required, but this is often not possible as $D$ increases with each iteration. However, we ensure that the following bound is valid:
\begin{equation}
	\label{eq:smallproperty}
	0 \leq M \leq D+1,
\end{equation}
such that we add at most $D+1$ nodes to the quadrature rule, consequently obtaining a rule of at most $N+D+2$ nodes.

The first step in obtaining such a quadrature rule is to consider the following quadrature rule:
\begin{equation}
\label{eq:firstquadraturerule}
\begin{aligned}
	X_{N+K+1} &= \{ \x_0, \dots, \x_N, \y_0, \dots, \y_K \}, \\
	W_{N+K+1} &= \Bigl\{ \underbrace{0, \vphantom{\frac{1}{K+1}}\dots, 0}_{N+1}, \underbrace{\frac{1}{K+1}, \dots, \frac{1}{K+1}}_{K+1} \Bigr\}.
\end{aligned}
\end{equation}
It is trivial to see that this quadrature rule is such that \eqref{eq:keyproperty} holds, but it is not such that \eqref{eq:smallproperty} holds. The idea is to iteratively remove nodes from this quadrature rule keeping \eqref{eq:keyproperty} as an invariant, which can be done until \eqref{eq:smallproperty} holds. All constructed intermediate quadrature rules have positive weights.

To start off, notice that \eqref{eq:firstquadraturerule} in conjunction with \eqref{eq:keyproperty} can be rewritten as the following linear system:
\begin{equation}
	\begin{pmatrix}
		\varphi_0(\x_0) & \cdots & \varphi_0(\x_N) & \varphi_0(\y_0) & \cdots & \varphi_0(\y_K) \\
		\vdots & \ddots & \vdots & \vdots & \ddots & \vdots \\
		\varphi_D(\x_0) & \cdots & \varphi_D(\x_N) & \varphi_D(\y_0) & \cdots & \varphi_D(\y_K)
	\end{pmatrix}
	\begin{pmatrix}
		w_0 \\
		\vdots \\
		w_N \\
		w_{N+1} \\
		\vdots \\
		w_{N+K+1}
	\end{pmatrix}
	=
	\begin{pmatrix}
		\mu_0 \\
		\vdots \\
		\mu_D
	\end{pmatrix},
\end{equation}
or denoted equivalently by:
\begin{equation}
	V_D(X_{N+K+1}) \mathbf{w} = \mub.
\end{equation}
Here $\mu_j$ are the moments determined by means of sampling (``sample moments''): $\mu_j = 1/(K+1) \sum_{k=0}^K \varphi_j(\y_k)$. The Vandermonde matrix of this linear system (see Section~\ref{subsubsec:degreepos}) is a $(D+1) \times (N+K+2)$-matrix, so it has $J \coloneqq (N+K+2)-(D+1)$ linearly independent null vectors $\mathbf{c}_1, \dots, \mathbf{c}_J$, with
\begin{equation}
	V_D(X_{N+K+1}) \mathbf{c}_j = \mathbf{0}, \text{ for $j = 1, \dots, J$}.
\end{equation}
Hence it holds that
\begin{equation}
	V_D(X_{N+K+1}) (\mathbf{w} - \mathbf{c}) = \mub, \text{ for any $\mathbf{c} \in \linspan\{\mathbf{c}_1, \dots, \mathbf{c}_J\}$}.
\end{equation}
These null vectors can be used to remove $J$ nodes from the quadrature rule $X_{N+K+1}$, $W_{N+K+1}$~\cite{Bos2016b}, which is done by constructing a $\mathbf{c} \in \linspan\{\mathbf{c}_1, \dots, \mathbf{c}_J\}$ such that:
\begin{subequations}
\begin{align}
	\mathbf{w} - \mathbf{c} &\geq 0, \label{eq:condition1} \\
	w_{k_j} - c_{k_j} &= 0, \text{ for $J$ distinct indices $k_1, \dots, k_J$}. \label{eq:condition2}
\end{align}
\end{subequations}
Condition \eqref{eq:condition1} describes that all weights should be non-negative and condition \eqref{eq:condition2} describes that there should be $J$ weights that equal zero. If such a vector $\mathbf{c}$ is constructed, the nodes $k_1, \dots, k_J$ can be removed from the quadrature rule without loss of accuracy. Then the weights can straightforwardly be determined by $\mathbf{w} \gets \mathbf{w} - \mathbf{c}$.

The algorithm to determine such a $\mathbf{c}$ is initiated with $\mathbf{c} = \mathbf{0}$ and $C = \{\mathbf{c}_1, \dots, \mathbf{c}_J\}$. At each iteration, $\mathbf{c}$ is updated such that firstly \eqref{eq:condition1} is still valid and secondly such that one more weight is equal to zero. Therefore, consider $\mathbf{c}_1$, an element (and without loss of generality, the first element) of $C$. Since $\mathbf{c}_1$ is a null vector, it holds that $\mathbf{c} + \alpha \mathbf{c}_1$ is a null vector for any $\alpha \in \mathbb{R}$. The goal is to determine $\alpha$ such that conditions \eqref{eq:condition1} and \eqref{eq:condition2} hold. Condition \eqref{eq:condition1} translates to
\begin{equation}
	\mathbf{w} - \mathbf{c} - \alpha \mathbf{c}_1 \geq 0,
\end{equation}
which is equivalent to
\begin{gather}
	 \alpha_\textrm{min} \leq \alpha \leq \alpha_\textrm{max}, \text{ with} \\
	 \alpha_\textrm{min} = \max\left( \frac{w_k - c_k}{c_{1,k}} \mid c_{1,k} < 0 \right) \text{ and } \alpha_\textrm{max} = \min\left( \frac{w_k - c_k}{c_{1,k}} \mid c_{1,k} > 0 \right).
\end{gather}
Here $w_k$, $c_k$, and $c_{1,k}$ denote the elements of $\mathbf{w}$, $\mathbf{c}$, and $\mathbf{c}_1$, respectively. Notice that $\alpha_\textrm{min} < 0 < \alpha_\textrm{max}$, since $c_{1,k} < 0$ for $\alpha_\textrm{min}$ and $c_{1,k} > 0$ for $\alpha_\textrm{max}$. Hence $\alpha$ is well-defined in this way. By choosing either $\alpha = \alpha_\textrm{min}$ or $\alpha = \alpha_\mathrm{min}$, there exists one $k_0$ such that
\begin{equation}
	w_{k_0} - c_{k_0} - \alpha c_{1,k_0} = 0.
\end{equation}
Therefore, $\mathbf{c}$ is updated such that $\mathbf{c} \gets \mathbf{c} + \alpha \mathbf{c}_1$ with $\alpha = \alpha_\textrm{min}$ or $\alpha = \alpha_\textrm{max}$.To ensure that the $k_0$-th weight remains zero in subsequent iterations, all other null vectors in $C$ are updated as follows:
\begin{equation}
	\mathbf{c}_j = \mathbf{c}_j - \frac{c_{j,k_0}}{c_{1,k_0}} \mathbf{c}_1, \text{ for all $j > 1$}.
\end{equation}
This is essentially a Gaussian reduction with pivot $c_{1,k_0}$ and ensures that $c_{j,k_0} = 0$ for all $j > 1$, which yields that $w_{k_0} - c_{k_0} = 0$ in all subsequent iterations (independently of the chosen $\alpha$ in that iteration). Finally $\mathbf{c}_1$ is removed from $C$ and the process is repeated until $C$ is empty. If $C$ is empty, a $\mathbf{c}$ is obtained that satisfies both conditions \eqref{eq:condition1} and \eqref{eq:condition2}.

By determining such a $\mathbf{c}$, the vector $\mathbf{w} - \mathbf{c}$ has $J$ coefficients that are equal to 0. We do not want to remove all $J$ nodes and weights, since the first $N+1$ weights belong to nodes where the costly computational model has been evaluated. A straightforward solution to this is to construct the following quadrature rule:
\begin{align}
	X_{N+M} &= \{ \x_k \mid k \leq N \text{ or } w_k - c_k \neq 0 \}, \\
	W_{N+M} &= \{ w_k - c_k \mid k \leq N \text{ or } w_k - c_k \neq 0 \}.
\end{align}
The vector $\mathbf{w} - \mathbf{c}$ has at least $J$ zero elements. So it has at most $D+1$ nonzero elements. Hence the quadrature rule based on $X_{N+M}$ and $W_{N+M}$ has at most $(N+1) + (D+1)$ nodes. Notice that therefore this is a quadrature rule that satisfies \eqref{eq:keyproperty} and \eqref{eq:smallproperty}.

It is not necessary to compute the full null space of the $(D+1) \times (N+K+2)$-matrix. Each iteration a null vector of the full matrix can be obtained by computing a null vector of a $(D+1) \times (D+2)$ submatrix using $D+2$ columns that represent nodes that are still present in the quadrature rule, and subsequently padding the obtained vector with elements equal to 0.

There is freedom in the value of $\alpha$, i.e.\ both $\alpha = \alpha_\textrm{max}$ and $\alpha = \alpha_\textrm{min}$ yield a valid quadrature rule with positive weights. We use this freedom to determine for each null vector $\mathbf{c}_k$ the $\alpha$ that yields the smallest quadrature rule (i.e.\ $w_k - c_k = 0$ for as many $k > N$ as possible).

\subsection{Construction of the proposal distribution}
\label{subsec:interpolation}
By using the quadrature rule constructed using the techniques from the previous section, a Bayesian prediction can be inferred by evaluating the model at the quadrature rule nodes and assessing the accuracy using the nesting of the quadrature rules. By using Bayes' rule, evaluations of the (unscaled) posterior are readily available, that is $\rho(\x_k)$ for $k = 0, \dots, N+M$. The goal is to use these evaluations to construct an approximation of the posterior.

Constructing an approximation of a distribution is non-trivial and there are some specific requirements on the interpolation algorithm used in this work. These requirements are further discussed in Section~\ref{subsubsec:prerequisites}. In this article, all results are obtained using nearest neighbor interpolation, which is explained in Section~\ref{subsubsec:nni}. There exist various alternative approaches, with varying advantages and disadvantages that warrant our choice for nearest neighbor interpolation. These alternatives are reviewed in Section~\ref{subsubsec:alternatives}.

\subsubsection{Prerequisites of the interpolant}
\label{subsubsec:prerequisites}
The goal is to construct an approximate posterior $\rho_{N+M}$ such that the available $N+M+1$ point evaluations of $\rho$ are taken into account. These point evaluations of the posterior are exact, so in this work we construct $\rho_{N+M}$ such that it is interpolatory, i.e.\
\begin{equation}
	\rho_{N+M}(\x_k) = \rho(\x_k), \text{ for all $k = 0, \dots, N+M$}.
\end{equation}
Obviously, the interpolant should be constructed such that $\rho_{N+M}(\x) \approx \rho(\x)$ for all $\x \neq \x_k$ (for any $k$). It is not necessary that $\rho_{N+M}$ is interpolatory for the construction of the quadrature rule, but it is assumed to be the case in the analysis discussed in Section~\ref{sec:erranalysis}.

It is known that $\rho$ is a probability density function, possibly up to a finite constant. Therefore the interpolant should be such that $\rho_{N+M}(\x) \geq 0$ for all $\x$. Moreover, since a quadrature rule is constructed using this distribution, all moments must be finite, i.e.
\begin{equation}
	\int_\Omega \varphi(\x) \, \rho_{N+M}(\x) \dd \x < \infty, \text{ for all polynomials $\varphi$}.
\end{equation}
These two conditions ensure that $\rho_{N+M}$ is again a probability density function (possibly up to a constant).

We emphasize that in general it is \emph{not} preferable that the moments of $\rho_{N+M}$ are equal to the moments of the quadrature rule, i.e.\ in general it should hold that
\begin{equation}
	\I^{\rho_{N+M}} \varphi_j = \int_\Omega \varphi_j(\x) \, \rho_{N+M}(\x) \dd \x \neq \int_\Omega \varphi_j(\x) \, \rho_N(\x) \dd \x = \sum_{k=0}^{N+M} w_k \, \varphi_j(\x_k) = \A^{\rho_N}_{N+M} \varphi_j, \text{ for $j = 0, \dots, D_i$}.
\end{equation}
Here, we use the notation $\I^{\rho_{N+M}}$ and $\A^{\rho_N}_{N+M}$ as defined in Section~\ref{subsec:quadrules}, i.e.
\begin{align}
	\I^{\rho_{N+M}} u &= \int_\Omega u(\x) \, \rho_{N+M}(\x) \dd \x, \\
	\A^{\rho_N}_{N+M} u &= \sum_{k=0}^{N+M} w_k \, u(\x_k).
\end{align}
The moments estimated by the quadrature rule are generally not equal to the exact moments of the true distribution $\rho$. Hence forcing that the moments of $\rho_{N+M}$ are equal to the moments of the (current) quadrature rule does not lead to a converging sequences of distributions (as all moments would be equal to the approximation of the first quadrature rule).

\subsubsection{Nearest neighbor interpolation}
\label{subsubsec:nni}
In this work, $\rho_{N+M}$ is constructed by means of nearest neighbor interpolation of the likelihood and multiplying the obtained interpolant with the prior. Hence for any $\x \in \Omega$ the interpolant $\rho_{N+M}(\x)$ is defined as follows:
\begin{equation}
	\rho_{N+M}(\x) = q(\z \mid \x_k) \, q(\x), \text{ with $\x_k = \argmin_{\x_j \in X_{N+M}} \| \x - \x_j \|_2$}, \text{ with } X_{N+M} = \{\x_0, \dots, \x_{N+M}\}.
\end{equation}
The advantages of this approach are among others that the obtained interpolant is always non-negative and globally bounded, i.e.\ for all $\x$ it holds that
\begin{equation}
	0 \leq \rho_{N+M}(\x) \leq \left(\sup_{\x \in \Omega} q(\z \mid \x)\right) \left(\sup_{\x \in \Omega} q(\x)\right).
\end{equation}
If an improper prior on an unbounded domain is considered, it might be the case that $\rho_{N+M}$ is not a well-defined distribution. This can be alleviated by restricting the nearest neighbor interpolant to a compact set that is larger than the convex hull of the quadrature rule nodes.

Moreover, for the purposes of this article it can be implemented very efficiently, as the number of interpolation points is relatively small and the exact structure of the interpolant is not of importance: sampling from $\rho_{N+M}$ can simply be done by means of acceptance rejection sampling using the prior as proposal distribution (which requires that the prior is not improper).

In particular, for the nearest neighbor interpolant it is not necessary to construct the Voronoi diagram of the nodes, but only to implement a routine that can find the closest node for any given input sample. This can be done very fast and efficiently for large numbers of samples. Moreover, searching for the nearest node has only a weak dependence of the dimension on the space (contrary to constructing the Voronoi diagram, which is very costly in higher dimensional spaces), which allows to use this approach in high dimensions.

\begin{figure}
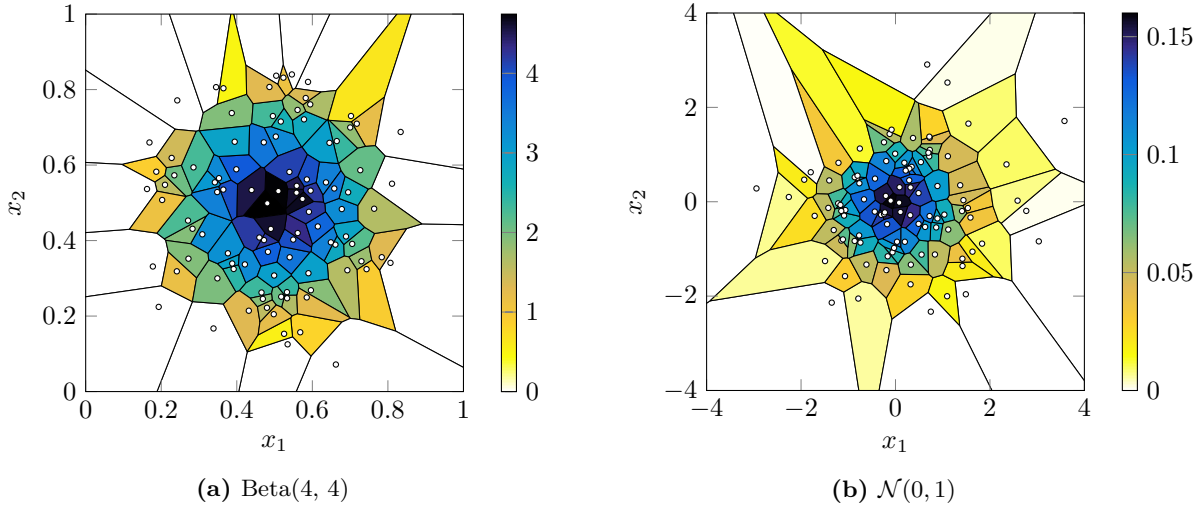

	\begin{minipage}{.5\textwidth}
		\centering
		\includepgf{.8\textwidth}{.8\textwidth}{nn-beta44.tikz}
		\subcaption{Beta(4, 4)}
	\end{minipage}%
	\begin{minipage}{.5\textwidth}
		\centering
		\includepgf{.8\textwidth}{.8\textwidth}{nn-normal.tikz}
		\subcaption{$\mathcal{N}(0, 1)$}
	\end{minipage}

	\caption{Nearest neighbor interpolations of random samples based on two distributions. In both cases, $x_1$ and $x_2$ are identically and independently distributed.}
	\label{fig:nn}
\end{figure}

Nearest neighbor interpolation is being used often in statistics~\cite{Biau2015}, especially with a focus on classification problems. It belongs to the broader class of scattered data interpolation methods~\cite{Franke1982} and is arguably one of the most straightforward methods to use in this regard. Moreover, it is a local interpolation method, whereas the quadrature rule is a global integration method, which naturally balances exploitation of locality with exploration of the space without any free parameters that have to be tuned.

As mentioned briefly before, a nearest neighbor interpolant can also be applied directly to the integral for the construction of quadrature rules~\mbox{\cite{Du1999,Guessab2008,Luo2009,Butler2017}}. Many of these approaches focus on explicitly constructing the nodal set such that the nodes are the centroids (centers of mass) of the cells of the Voronoi diagram. These approaches cannot be extended straightforwardly to the setting to this article, since it is often assumed that the distribution is computationally accessible. Moreover, it is usually necessary to construct the Voronoi diagram, which hinders the applicability to higher dimensional spaces.

The idea of using a nearest neighbor interpolant to construct a surrogate has been explored successfully before and various improvements can be considered to enhance the accuracy of the interpolant~\mbox{\cite{Rushdi2017,Ebeida2016}}. The advantages of constructing a surrogate in this way are among others that the surrogate exploits locality, as mentioned above. The main difference between the approaches used in \mbox{\cite{Rushdi2017,Ebeida2016}} and the approach considered in this article is that in this work a surrogate of a probability distribution is generated, which adds several constraints to the approach (see Section~\ref{subsubsec:prerequisites}).

Nearest neighbor interpolation has been illustrated in Figure~\ref{fig:algorithm} (see page~\pageref{fig:algorithm}).Two other examples of interpolated probability density functions are depicted in Figure~\ref{fig:nn}. For sake of illustration, the interpolation nodes of Figure~\ref{fig:nn} are based on random samples.

\subsubsection{Alternative distribution approximation methods}
\label{subsubsec:alternatives}
It is non-trivial to consider more accurate interpolation methods in this framework, as assuring that the obtained interpolant is a PDF is difficult. For example, many global methods such as Lagrange interpolation and Gaussian process regression do not yield an interpolant that is non-negative. Extensions to enforce positivity exist, for example for Lagrange interpolation~\cite{Grzelak2018,Fritsch1980}. Often these approaches only work in univariate cases. Gaussian process regression can possibly be used by using a strictly positive prior that is not improper and enforcing finiteness of integrals over the obtained approximation~\cite{Rasmussen2006}, but this is rather involved and conflicts with the convenient Bayesian framework behind Gaussian processes.

Another well-known method to construct a suitable approximation is the kernel density estimate, or more general, the family of methods based on radial basis functions~\cite{Buhmann2000}. However, in the cases in this article explicit values of the probability density function are known, which cannot be incorporated in a straightforward fashion in these methods: enforcing specific values of a function in a radial basis function setting does not necessarily yield a positive interpolant.

A promising method that balances the accuracy of polynomial interpolation and piecewise behavior of nearest neighbor interpolation is the Simplex Stochastic collocation method~\cite{Witteveen2010,Witteveen2012,Witteveen2013,Edeling2016}, where piecewise high order polynomial interpolation is combined with a Delaunay triangulation of the nodal set. This method converges exponentially for sufficiently smooth functions and is non-oscillatory (i.e.\ it is bounded from below and above by the minimum and maximum of the model), so it is very suitable for the interpolation of a PDF. However, it requires the explicit construction of the interpolant with the underlying Delaunay triangulation, which severely limits its applicability to multi-dimensional spaces.

Concluding, it is nontrivial to consider a more efficient alternative approach to nearest neighbor interpolation that yields a positive interpolant with similar efficiency as nearest neighbor interpolation. It is out of the scope of this work to derive such a new interpolation approach. Moreover, nearest neighbor interpolation is sufficiently accurate for the examples considered in this work, as demonstrated in Section~\ref{sec:numerics}.

\section{Error analysis}
\label{sec:erranalysis}
The accuracy of the proposed method is studied by assessing the convergence of the sequence of quadrature rules. The overall error we are interested in is the integration error, defined for a given continuous function $f: \Omega \to \mathbb{R}$ as follows:
\begin{equation}
	\label{eq:integrationerror}
	e_N(f) \coloneqq | \I^\rho f - \A^{\rho_N}_N f |.
\end{equation}
Here the operator $\A^{\rho_N}_N$ has been constructed using the methods outlined in Section~\ref{sec:quadraturerule}, i.e.\ for known $D_i$ we have
\begin{equation}
	\A^{\rho_N}_N \varphi_j = \sum_{k=0}^N w_k \, \varphi_j(\x_k) = \int_\Omega \varphi_j(\x) \, \rho_N(\x) \dd \x = \I^{\rho_N} \varphi_j, \text{ for $j = 0, \dots, D_i$},
\end{equation}
where we neglect the effect of using samples for the construction of the quadrature rule.

It can be demonstrated that $e_N \to 0$ if $N \to \infty$, provided that $f$ is sufficiently smooth and $D_i \to \infty$ for $i \to \infty$. The details of this are discussed in Section~\ref{subsec:convergence}.

A central question is whether this error decays faster than the integration error of a conventional quadrature rule constructed using the prior. As constructing such a quadrature rule does not require interpolation or sampling, there are significantly fewer sources of error. This topic can be addressed by using principles from importance sampling, which is outlined in Section~\ref{subsec:errimportance}.

The error analysis based on importance sampling neglects any error that occurs due to the usage of samples of $\rho_{N+M}$ instead of the full distribution. We do not consider this a severe limitation, since sampling from $\rho_{N+M}$ does not require costly model evaluations and therefore many samples can be used to minimize this error. Nonetheless, this error will be briefly considered in Section~\ref{subsec:errsampling}.

\subsection{Decay of the integration error}
\label{subsec:convergence}
Notice that the integration error $e_N(f)$ can be decomposed into three components: an interpolation error, a sampling error, and a quadrature rule error:
\begin{align}
	e_N(f) = | \I^\rho f - \A^{\rho_N}_N f | &= | \I^\rho f - \I^{\rho_N} f + \I^{\rho_N} f - \A^{\rho_N}_N f | \\
	&\leq | \I^\rho f - \I^{\rho_N} f | + | \I^{\rho_N} f - \A^{\rho_N}_N f | \\
	&= | \I^\rho f - \I^{\rho_N} f | + | \I^{\rho_N} f - \I^{\rho_N}_K f + \I^{\rho_N}_K f - \A^{\rho_N}_N f | \\
	&\leq \underbrace{| \I^\rho f - \I^{\rho_N} f |}_{\mathclap{\text{Interpolation error $I_N$}}} ~+~ \underbrace{| \I^{\rho_N} f - \I^{\rho_N}_K f |}_{\mathclap{\text{Sampling error $S_K$}}} ~+~ \underbrace{| \I^{\rho_N}_K f - \A^{\rho_N}_N f |}_{\mathclap{\text{Quadrature error $Q_N$}}} \\
	&= I_N + S_K + Q_N.
	\label{eq:basicconvergence}
\end{align}
\begin{equation}
	\I^{\rho_N}_K f \coloneqq \frac{1}{K+1} \sum_{k=0}^K f(\y_k).
\end{equation}
By introducing this operator, the quadrature rule error $Q_N$ can be bounded using the Lebesgue inequality from \eqref{eq:lebesgue} as discussed in Section~\ref{subsubsec:degreepos}:
\begin{equation}
	Q_N \leq \left(1 + \sum_{k=0}^N |w_k| \right) \inf_{\varphi \in \Phi_{D_i}} \| f - \varphi \|_\infty.
\end{equation}
Since all quadrature rules constructed in this work have positive weights, this error decays if $D_i \to \infty$ (with $i \to \infty$) for sufficiently smooth $f$.

The interpolation error $I_N$ depends on the specific procedure used to construct the interpolant. It holds that $I_N \to 0$ for $N \to \infty$ if $\| \rho_N - \rho \|_\infty \to 0$. The latter can be demonstrated straightforwardly for nearest neighbor interpolants by demonstrating that the quadrature rule nodes distribute across the space, i.e.\ the mutual distance between two closest nodes decays to zero. A nearest neighbor interpolant yields a converging interpolant in that case.

To demonstrate that quadrature rule nodes distribute across the space, let $\x_0, \dots, \x_k, \dots$ be a sequence of nodes such that $\{\x_0, \dots, \x_N\}$ form the nodes of a quadrature rule with positive weights with respect to a distribution $\rho$ for all $N$. Moreover assume that $\x_k \in \Omega$ for all $k \in \mathbb{N}$. Suppose $S \subset \Omega$ is such that there is no $k \in \mathbb{N}$ such that $\x_k \in S$. To demonstrate that the nodes distribute across the space, we will show that $S$ as considered here is a null set, i.e.
\begin{equation}
	\int_S \rho(\x) \dd \x = 0,
\end{equation}
which implies that the mutual distance between the quadrature rule nodes vanishes. To demonstrate that $S$ is a null set, let $u_S \in C^\infty(\Omega)$ be a smooth function such that $u_S(\x) = 0$ for all $\x \not\in S$ and such that
\begin{equation}
	\int_\Omega u_S(\x) \, \rho(\x) \dd \x > 0.
\end{equation}
Such a $u_S$ only exists if $S$ is not a null set. A univariate example of such a function is the bump function, e.g.\ if $S = [-1, 1]$ then
\begin{equation}
	u_S(x) = \begin{cases}
		\exp\left( -\frac{1}{1-x^2} \right) &\text{if $|x| < 1$}, \\
		0 &\text{otherwise}.
	\end{cases}
\end{equation}
Similarly, multivariate bump functions can be constructed by products of univariate bump functions defined on a hypercube in $S$. Due to the smoothness of $u_S$, there exists a sequence of polynomials $\varphi_N$ such that $\| u_S - \varphi_N \|_\infty \to 0$ for $N \to \infty$. Hence \emph{any} quadrature rule operator $\A^\rho_N$ with positive weights yields a converging estimate of $\I^\rho u_S$, i.e.
\begin{equation}
	\A^\rho_N u_S = \sum_{k=0}^N w_k \, u_S(\x_k) \to \int_\Omega u_S(\x) \, \rho(\x) \dd \x, \text{ for $N \to \infty$}.
\end{equation}
We assumed that the sequence of nodes is such that there is no $\x_k \in S$ for any $k \in \mathbb{N}$. Since the quadrature rule converges, this yields that the integral of $u_S$ must vanish:
\begin{equation}
	\int_\Omega u_S(\x) \, \rho(\x) \dd \x = \int_S u_S(\x) \, \rho(\x) \dd \x = 0.
\end{equation}
This holds for any smooth function defined in this way, which is only possible if
\begin{align}
	\int_S \rho(\x) \dd \x &= 0.
\end{align}
Notice that positivity of the weights is essential here, as otherwise it is not guaranteed that the quadrature rule approximation of any smooth function converges.

A similar derivation can be formulated for the nodes deduced with the methodology proposed in this article. To see this, notice that for $u_S$, as introduced above, it holds that
\begin{equation}
	\left| \sum_{k=0}^N w_k \, u_S(\x_k) - \int_\Omega u_S(\x) \, \rho_N(\x) \dd \x \right| \leq \inf_{\varphi \in \Phi_{D_i}} \| u_S - \varphi \|_\infty,
\end{equation}
which simplifies to (recall that $u_S(\x_k) = 0$ for all $k$)
\begin{equation}
	\label{eq:lebesguerhon}
	\left| \int_S u_S(\x) \, \rho_N(\x) \dd \x \right| \leq \inf_{\varphi \in \Phi_{D_i}} \| u_S - \varphi \|_\infty.
\end{equation}
These inequalities hold for any $N$. Recall that $\rho(\x) > 0$ for all $\x$, and therefore $\rho_N(\x) > 0$ for all $\x$ and $N$. Hence, if $\int_S \rho(\x) \dd \x > 0$, the left hand side of \eqref{eq:lebesguerhon} is bounded from below by a positive constant that is independent of $N$. Following the same argument as above, we therefore conclude that $S$ is a null set.

The sampling error $S_K$ is independent from $N$, and can therefore be reduced without any additional costly model evaluations by considering more samples from $\rho$. The exact convergence rate depends on the constructed sequence of samples, which will be further discussed in Section~\ref{subsec:errsampling}.

Concluding, if $i \to \infty$ and $N \to \infty$, it holds that $e_N(f) \to S_K$ provided that $\inf_{\varphi \in \Phi_{D_i}} \| \varphi - f \|_\infty \to 0$. Moreover, if also $K \to \infty$, it holds that $e_N(f) \to 0$. This demonstrates the convergence of our method for sufficiently smooth functions, but this does not demonstrate that our method converges faster than constructing a conventional quadrature rule with respect to the prior, since such a rule does not have an interpolation error. This is the main topic of Section~\ref{subsec:errimportance}.

\subsection{Importance sampling}
\label{subsec:errimportance}
The proposed quadrature rule of this article has three sources of error: an interpolation error $I_N$, a sampling error $S_K$, and a quadrature rule error $Q_N$. A quadrature rule constructed using the prior (e.g.\ by means of a sparse grid) only has the quadrature rule error, i.e.\ \eqref{eq:basicconvergence} simplifies to $e_N(f) = Q_N$.

The main difference is that a rule that only uses the prior requires a different integrand, i.e.\ $\I^\rho f$ is approximated by means of $\A^{q(.)}_N (q(\z \mid \x) f(\x))$, whereas the quadrature rule proposed in this work uses $f(\x)$ as integrand. In this section we derive in which cases the proposed quadrature rule needs fewer model evaluations than a conventional quadrature rule to reach a similar accuracy. Without loss of generality, we assume that $\rho$ and $\rho_N$ are both probability density functions, which yields that $\| \mathcal{I}^\rho \|_\infty = \| \mathcal{I}^{\rho_N} \|_\infty = 1$.

The analysis is based on the main principle of importance sampling, a technique that is used to estimate integrals with respect to a target distribution if samples from a proposal distribution are known. Even though the main interest is not in drawing samples from the posterior, the following key idea of importance sampling can be readily used: let $\psi(\x)$, the so-called proposal distribution, be a probability density function with the same support as $\rho$. Then
\begin{equation}
	\I^\rho f = \int_\Omega f(\x) \, \rho(\x) \dd \x = \int_\Omega f(\x) \frac{\rho(\x)}{\psi(\x)} \, \psi(\x) \dd \x = \I^\psi \left( f \frac{\rho}{\psi} \right).
\end{equation}
Hence it is not necessary to calculate a possibly expensive integral with respect to the distribution $\rho$, but to calculate one with respect to the proposal distribution $\psi$.

\begin{figure}
	\centering
	\includegraphics[width=.4\textwidth]{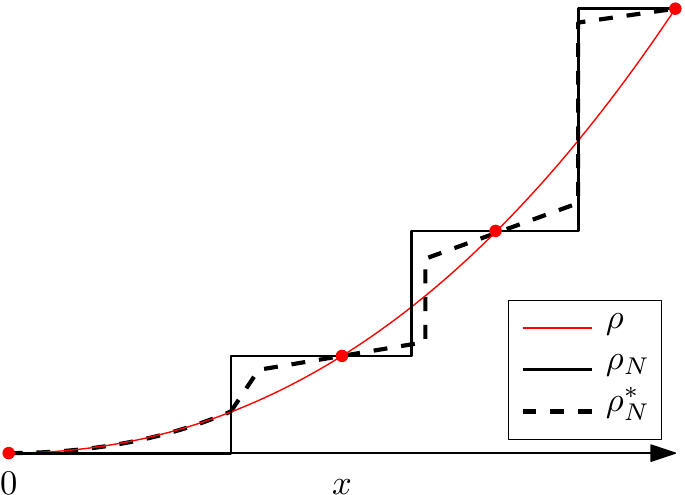}
	\caption{A sketch of the motivation why $\rho^*_N$ is introduced. In this particular case, $\rho(x) / \rho_N(x)$ is unbounded for $x \to 0$. Therefore, a $\rho^*_N$ is introduced that yields a bounded fraction near $x = 0$. Small corrections are necessary to enforce that the raw moments do not change.}
	\label{fig:rho-star}
\end{figure}

It is a natural idea to use the proposal distribution $\psi = \rho_N$ for the analysis, but this can severely limit the applicability of the Lebesgue inequality from \eqref{eq:lebesgue}, as $\| \rho / \rho_N \|_\infty$ can be large even if $\| \rho - \rho_N \|_\infty$ is small. In many cases, $\rho / \rho_N$ can be globally bounded. For example, if a uniform prior and a Gaussian likelihood is considered, $\rho_N$ can be used as proposal distribution. An example where this is not the case is if $\rho$ is a Beta distribution and a node is placed on the boundary of $\Omega$. In that case $\rho / \rho_N$ is unbounded on the boundary.

To alleviate this, we choose a slightly different approach, and use a proposal distribution $\rho^*_N$ such that it is interpolatory and integrates the same moments as $\rho_N$, i.e.
\begin{subequations}
\begin{equation}
	\label{eq:rhostar1}
	\rho^*_N(\x_k) = \rho_N(\x_k) = \rho(\x_k), \text{ for $k = 0, \dots, N$},
\end{equation}
and
\begin{equation}
	\label{eq:rhostar2}
	\I^{\rho^*_N} \varphi_j = \I^{\rho_N} \varphi_j, \text{ for $j = 0, \dots, D_i$}.
\end{equation}
\end{subequations}
These two conditions define a set of possible proposal distributions, which we will call $R$ in this section, i.e.\ $\rho^*_N \in R$ if it solves \eqref{eq:rhostar1} and \eqref{eq:rhostar2}. See Figure~\ref{fig:rho-star} for a sketch how such a $\rho^*_N$ can be interpreted. In this example, $\rho^*_N$ is discontinuous and roughly following $\rho_N$, but this is not necessary for the analysis. The only conditions on $\rho^*_N$ are \eqref{eq:rhostar1}, \eqref{eq:rhostar2}, and $\| \rho / \rho^*_N \|_\infty < \infty$. It is actually not necessary to construct a $\rho^*_N$ to use the method proposed in this work: only samples from $\rho_N$ are used. The distribution $\rho^*_N$ is only a tool in the analysis discussed in this section. Notice that $\rho_N \in R$, so $R$ is a well-defined non-empty set. However, as explained above it is not desirable to use $\rho^*_N = \rho_N$ in the subsequent analysis as that leads to a potentially large $\| \rho / \rho^*_N \|_\infty$.

The two conditions stated above ensure firstly that the moments of $\rho_N$ are equal to the moments of $\rho^*_N$. Secondly, integrating $\rho_N$, $\rho$, or $\rho^*_N$ using a quadrature rule with nodes $\x_0, \dots, \x_N$ yields the same result. The latter is especially useful in importance sampling, since this implies that for all $k = 0, \dots, N$:
\begin{equation}
	1 = \frac{\rho(\x_k)}{\rho_N(\x_k)} = \frac{\rho(\x_k)}{\rho^*_N(\x_k)}.
\end{equation}
This property will be used extensively in the remainder of this section. To this end, choose a $\rho^*_N \in R$, i.e.\ one that satisfies \eqref{eq:rhostar1} and \eqref{eq:rhostar2}, as proposal distribution. Then
\begin{equation}
	\label{eq:ingredients1}
	\I^{\rho} f = \I^{\rho_N} \left( f \frac{\rho}{\rho_N} \right) = \I^{\rho^*_N} \left( f \frac{\rho}{\rho^*_N} \right).
\end{equation}
A similar argument can be applied to rewrite $\A_N^{\rho_N}$:
\begin{equation}
	\A_N^{\rho_N} f = \sum_{k=0}^N w_k f(\x_k) = \sum_{k=0}^N w_k \frac{\rho(\x_k)}{\rho_N(\x_k)} f(\x_k) = \A_N^{\rho_N} \left( f \frac{\rho}{\rho_N} \right),
\end{equation}
where it is being used that $\rho_N(\x_k) = \rho(\x_k)$ for all $k = 0, \dots, N$. Hence the following is obtained:
\begin{equation}
	\label{eq:ingredients2}
	\A_N^{\rho_N} f = \A_N^{\rho_N} \left( f \frac{\rho}{\rho_N} \right) = \A_N^{\rho_N} \left( f \frac{\rho}{\rho^*_N} \right).
\end{equation}

These expressions, i.e.\ the defining properties of $\rho^*_N$, as described by \eqref{eq:rhostar1} and \eqref{eq:rhostar2}, and the derived relations for the integral and quadrature rule operator, as described by \eqref{eq:ingredients1} and \eqref{eq:ingredients2}, form sufficient ingredients to bound $e_N(f)$. The derivation is similar to the derivation of the Lebesgue inequality (recall \eqref{eq:lebesgue}), but some bookkeeping is required to keep track of the correct proposal distribution:
\begin{align}
	e_N(f) &= | \I^\rho f - \A^{\rho_N}_N f | \\
	&= \left| \I^\rho f - \A^{\rho_N}_N \left( f \frac{\rho}{\rho^*_N} \right) \right| \\
	&= \left| \I^\rho f - \I^{\rho^*_N} \varphi + \I^{\rho^*_N} \varphi - \A^{\rho_N}_N \left( f \frac{\rho}{\rho^*_N} \right) \right|. \\
\intertext{Here, $\varphi \in \Phi_{D_i}$. Hence the quadrature rule is exact for $\varphi$ and the following is obtained:}
	e_N(f) &= \left| \I^{\rho^*_N} \left( f \frac{\rho}{\rho^*_N} \right) - \I^{\rho^*_N} \varphi + \A^{\rho_N}_N \varphi - \A^{\rho_N}_N \left( f \frac{\rho}{\rho^*_N} \right) \right|, \\
\intertext{where we use that $\I^{\rho^*_N} \varphi = \I^{\rho_N} \varphi = \A_N^{\rho_N} \varphi$ for $\varphi \in \Phi_{D_i}$, which follows from \eqref{eq:rhostar2}. Concluding, the following inequality is obtained:}
	e_N(f) &\leq \left| \I^{\rho^*_N} \left( f \frac{\rho}{\rho^*_N} \right) - \I^{\rho^*_N} \varphi \middle| + \middle| \A^{\rho_N}_N \varphi - \A^{\rho_N}_N \left( f \frac{\rho}{\rho^*_N} \right) \right| \\
	&\leq \left( \| \I^{\rho^*_N} \|_\infty + \| \A^{\rho_N}_N \|_\infty \right) \left\| f \frac{\rho}{\rho^*_N} - \varphi \right\|_\infty \\
	&= 2 \left\| f \frac{\rho}{\rho^*_N} - \varphi \right\|_\infty.
	\label{eq:ultimeineq}
\end{align}
Notice that in a Bayesian framework the ratio $\rho / \rho^*_N$ can be expressed independently from the prior:
\begin{equation}
	\frac{\rho(\x)}{\rho^*_N(\x)} = \frac{q(\z \mid \x) \, q(\x)}{q^*_N(\z \mid \x) \, q(\x)} = \frac{q(\z \mid \x)}{q^*_N(\z \mid \x)}.
\end{equation}
Here, with a little abuse of notation, $q^*_N(\z \mid \x)$ is the nearest neighbor interpolant of the likelihood. Similarly to $\rho^*_N$, it is not necessary to gain any knowledge about $q^*_N(\z \mid \x)$, since it is only used here to demonstrate independence from the prior.

By choosing $\rho^*_N \in R$ and $\varphi \in \Phi_{D_i}$ such that the obtained norm from \eqref{eq:ultimeineq} is minimal, an inequality similar to the Lebesgue inequality from \eqref{eq:lebesgue} is obtained:
\begin{equation}
	\label{eq:lebesgueimpr}
	e_N(f) = | \I^\rho f - \A^{\rho_N}_N f | \leq 2 \inf_{\rho^*_N} \inf_{\varphi \in \Phi_{D_i}} \left\| f \frac{\rho}{\rho^*_N} - \varphi \right\|_\infty = 2 \inf_{q^*_N} \inf_{\varphi \in \Phi_{D_i}} \left\| f \frac{q(\z \mid \cdot)}{q^*_N(\z \mid \cdot)} - \varphi \right\|_\infty.
\end{equation}
Here, $\rho^*_N(\x) = q^*_N(\z \mid \x) \, q(\x) \in R$ is according to \eqref{eq:rhostar1} and \eqref{eq:rhostar2}.

It is instructive to compare the obtained inequality with that of a conventional quadrature rule constructed with respect to the prior. If such a quadrature rule, say $\A^{q(\cdot)}_N$, is given, then by using the Lebesgue inequality its error can be bounded as follows (notice the difference with \eqref{eq:ultimeineq}):
\begin{equation}
	\label{eq:lebesgueconv}
	| \I^\rho f - \A^{q(\cdot)}_N \left( f \, q(\z \mid \cdot) \right) | = | \I^{q(\cdot)} \left( f \, q(\z \mid \cdot) \right) - \A_N^{q(\cdot)} \left( f \, q(\z \mid \cdot) \right) | \leq 2 \inf_{\varphi \in \Phi_{D_i}} \| f \, q(\z \mid \cdot) - \varphi \|_\infty,
\end{equation}
where we use that the prior is not improper (as otherwise $\| \I^{q(\cdot)} \|_\infty \neq 1$).

It is not generally true that the right hand side of \eqref{eq:lebesgueimpr}, i.e.\ the error of our proposed adaptive rule, is smaller than the right hand side of \eqref{eq:lebesgueconv}, i.e.\ the error of conventional rules with respect to the prior. For example, the latter is smaller if both $f$ and $q(\z \mid \cdot)$ are polynomials (which is rarely the case in Bayesian model calibration). In general, it is the case that a quadrature rule constructed using the iterative procedure of this article outperforms a quadrature rule constructed using the prior if $f$ can be approximated much more efficiently by polynomials than $f q(\z \mid \cdot)$. To see this, we assume that the approximation of $f$ converges \emph{strictly} faster than the approximation of $f q(\z \mid \cdot)$, so using little-o notation:
\begin{equation}
	\label{eq:fvsfq}
	\inf_{\varphi \in \Phi_{D_i}} \| f - \varphi \|_\infty = o\left(\inf_{\varphi \in \Phi_{D_i}} \| f \, q(\z \mid \cdot) - \varphi \|_\infty\right).
\end{equation}
To compare both quadrature rules, recall that $\rho_N$ is constructed such that $\rho_N \to \rho$. Then there exists a $\rho^*_N \in R$ such that
\begin{equation}
	\| \rho / \rho^*_N \|_\infty \to 1, \text{ for $N \to \infty$}.
\end{equation}
Hence there exist sequences $\{ a_N \}_N$ and $\{ b_N \}_N$ such that $a_N \leq \| \rho / \rho^*_N \|_\infty \leq b_N$ for all $N$ with $a_N \uparrow 1$ and $b_N \downarrow 1$ for $N \to \infty$\footnote{In other words, $a_N$ converges to $1$ from below and $b_N$ converges to $1$ from above for $N \to \infty$.}. Hence
\begin{align}
	\inf_{\varphi \in \Phi_{D_i}} \left\| f \frac{\rho}{\rho^*_N} - \varphi \right\|_\infty &\leq \max\left( \inf_{\varphi \in \Phi_{D_i}} \| a_N f - \varphi \|_\infty;~ \inf_{\varphi \in \Phi_{D_i}} \| b_N f - \varphi \|_\infty \right) \\
	&\leq \max\left( a_N \inf_{\varphi \in \Phi_{D_i}} \| f - \varphi \|_\infty;~ b_N \inf_{\varphi \in \Phi_{D_i}} \| f - \varphi \|_\infty \right).
\end{align}
So for all $\vareps > 0$ there exists an $N_1$ such that
\begin{equation}
	\inf_{\varphi \in \Phi_{D_i}} \left\| f \frac{\rho}{\rho_N^*} - \varphi \right\|_\infty \leq (1 + \vareps) \inf_{\varphi \in \Phi_{D_i}} \| f - \varphi \|_\infty, \text{ for all $N \geq N_1$}.
\end{equation}
Thus, for $N \geq N_1$ the following bound is obtained:
\begin{equation}
	\label{eq:forlargeN}
	e_N(f) \leq 2 (1+\vareps) \inf_{\varphi \in \Phi_{D_i}} \| f - \varphi \|_\infty.
\end{equation}
By combining this bound with \eqref{eq:fvsfq}, it follows that \eqref{eq:lebesgueimpr} is sharper than \eqref{eq:lebesgueconv}. In other words, we have shown that the decay of the error of the adaptive quadrature rule proposed in this work, denoted by $e_N(f)$ in \eqref{eq:forlargeN}, depends solely on whether $f$ can be approximated well by polynomials. Recall that whether $f$ can be approximated well using polynomials depends mostly on its smoothness properties, see Section~\ref{subsubsec:degreepos} for more discussion on this topic.

Concluding, if incorporating the likelihood in the integrand (obtaining $f q(\z \mid \cdot)$) yields an integrand that is difficult to approximate using polynomials, the proposed adaptive method has a sharper error bound than conventional rules. This is described by \eqref{eq:fvsfq} and is often the case if the likelihood is informative (i.e.\ it has small standard deviation) or if the model is irregular (which yields an even more irregular likelihood). For large $N$ the error of the adaptive quadrature rule proposed in this work behaves similar to that of a conventional quadrature rule constructed with respect to the \emph{posterior}, as described by \eqref{eq:forlargeN}, even though this conventional quadrature rule can usually not be determined without requiring large numbers of evaluations of the posterior.

\subsection{Sampling error}
\label{subsec:errsampling}
As an arbitrary sized sample set can be easily drawn from the distribution $\rho_N$ for any $N$, the sampling error at least decays with rate $\sqrt{K}$, i.e.\ for $K \to \infty$, the error converges to a Gaussian distribution as follows:
\begin{equation}
	| \I^{\rho_N} f - \I^{\rho_N}_K f | \sim \mathcal{N}(0, \sigma^2),
\end{equation}
with $\sigma = \sigma(f) / \sqrt{K}$. Here, $\sigma(f)$ is the standard deviation of $f$ with respect to $\rho_N$, i.e.
\begin{equation}
	\left(\sigma(f)\right)^2 = \I^{\rho_N} \bigl((f - \I^{\rho_N} f)^2\bigr).
\end{equation}
In this article, acceptance rejection sampling is used to draw large numbers of samples from $\rho_N$. Several improvements exist to construct sequences that convergence faster, for example quasi-Monte Carlo sequences~\cite{Caflisch1998,Niederreiter1992}. We emphasize that these approaches can be used straightforwardly in the framework explained in this article.

The number of samples $K$ can be chosen independently from the number of nodes $N$ under consideration, so the sampling error can be made arbitrarily small without any costly model evaluations. Choosing a larger sample set does however mean that it takes more computational effort to determine the implicit quadrature rule. In practice, the interpolation and quadrature rule error dominate both the error of the quadrature rule and the cost of the procedure (through model evaluations). Only if the model is analytic and the likelihood is not informative, the quadrature rule error and interpolation error decay rapidly enough to make the sampling error dominate.

\section{Numerical experiments}
\label{sec:numerics}
The properties of the method derived in this work are illustrated using various numerical examples. Two types of cases are discussed.

In Section~\ref{subsec:explicit} four different classes of analytical test cases are discussed to illustrate the key properties of our method. Firstly, in Section~\ref{subsubsec:domination} a univariate case is discussed, where the sampling and interpolation error as introduced in \eqref{eq:basicconvergence} can be observed well. Secondly, in Section~\ref{subsubsec:peakypeaky} three two-dimensional test functions are calibrated using numerically generated data. The test functions, based on Genz test functions~\cite{Genz1984}, are chosen such that the likelihood is very informative. Thirdly, in Section~\ref{subsubsec:genz} the six Genz test functions are calibrated using numerically generated data. The domain of these functions is five-dimensional and a comparison is made with conventional quadrature rules such as the tensor grid and the Smolyak sparse grid. Finally, it is instructive to assess the affect of varying dimension on the error of the quadrature rule. Therefore, in Section~\ref{subsubsec:genzdim} the calibration of the six Genz test functions is considered for varying dimension of $\Omega$.

The main motivation for this work is to apply Bayesian calibration to the closure coefficients of turbulence models. To demonstrate the applicability of our method to this case, we calibrate the closure coefficients of the Spalart--Allmaras one-equation model using measurement data from the flow over a transonic airfoil. This problem is further considered in Section~\ref{subsec:rae2822}.

\subsection{Analytical test cases}
\label{subsec:explicit}

\subsubsection{Effect of sampling and interpolation error}
\label{subsubsec:domination}
\begin{figure}
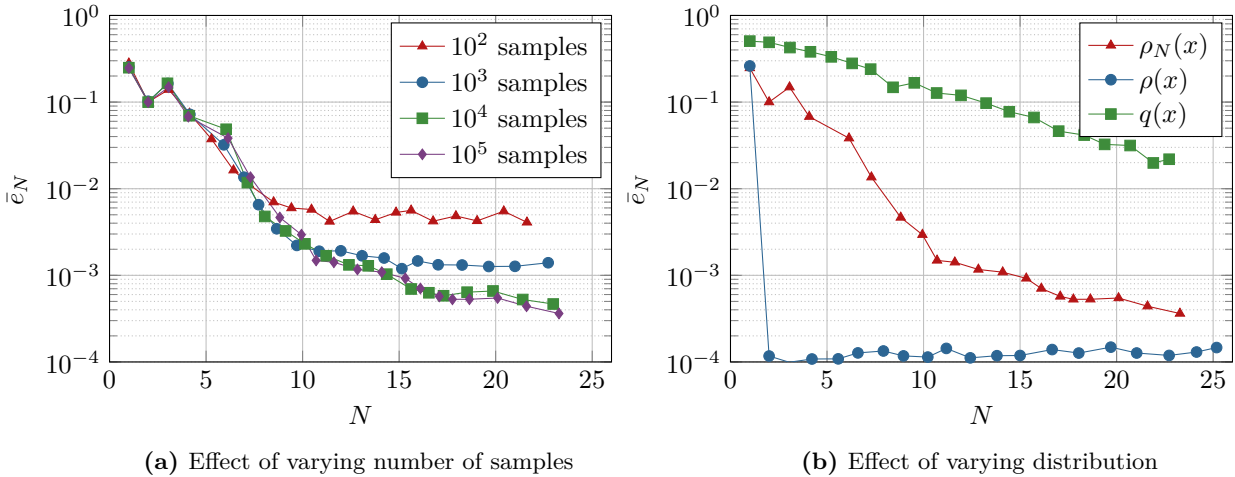

	\begin{minipage}[b]{.5\textwidth}
		\centering
		\includepgf{\textwidth}{.75\textwidth}{sampling-error.tikz}
		\subcaption{Effect of varying number of samples}
		\label{subfig:samplingerror}
	\end{minipage}%
	\begin{minipage}[b]{.5\textwidth}
		\centering
		\includepgf{\textwidth}{.75\textwidth}{interpolation-error.tikz}
		\subcaption{Effect of varying distribution}
		\label{subfig:interpolationerror}
	\end{minipage}

	\caption{Integration error as a function of the number of nodes. In the left figure, the error of an adaptive quadrature rule is depicted for various numbers of samples. In the right figure, the error of an implicit quadrature rule constructed using various distributions is depicted, for fixed $K = 10^5$.}
\end{figure}

Let $\Omega = [0, 1]$ and consider a uniformly distributed prior in conjunction with a $\mathrm{Beta}(40, 60)$-distributed likelihood, which yields:
\begin{align}
	\rho(x) &\propto x^{40}(1-x)^{60}.
\end{align}
To assess the effect of the \emph{sampling} error, the proposed adaptive quadrature rule is constructed using $K = 10^n$ samples with $n = 2, 3, 4, 5$. The obtained quadrature rule is used to approximate the mean of $\rho(x)$:
\begin{equation}
	\label{eq:estmean1}
	\int_\Omega \varphi_1(x) \, \rho(x) \dd x = \int_\Omega x \, \rho(x) \dd x \approx \sum_{k=0}^N w_k \, x_k = \sum_{k=0}^N w_k \, \varphi_1(x_k).
\end{equation}
The obtained estimate is compared to the true mean (which is $2/5$). The sampling error depends on the randomly generated samples, so the experiment is repeated 50 times and the obtained errors are averaged to illustrate general behavior of the error. The results are gathered in Figure~\ref{subfig:samplingerror}, where the averaged error is denoted by $\bar{e}_N$.

For small $N$, it is clear that the errors of all quadrature rules converge geometrically, which implies that the quadrature error or the interpolation error dominates. For large $N$, the error does not decay further due to the limited number of samples. The point where the error stops decaying depends on the number of samples used, i.e.\ a larger number of samples results in a smaller overall error.

To assess the effect of the \emph{interpolation} error, three different implicit quadrature rules are constructed. The first rule is constructed using the adaptive method proposed in this work, i.e.\ using $\rho_N(x)$ like in Figure~\ref{subfig:samplingerror}. This quadrature rule has a sampling error and interpolation error, but no quadrature error since we are comparing the first moment (see \eqref{eq:estmean1}). The second rule is constructed using the \emph{exact} distribution, i.e.\ $\rho(x)$, which results in a rule that only has a sampling error. The third rule is constructed using the uniform distribution, i.e.\ using $q(x)$. This rule has a sampling error and a quadrature error, as the likelihood needs to be incorporated in the integrand. All three quadrature rules are used to approximate the mean of the posterior. The rules constructed using $\rho_N(x)$ and $\rho(x)$ can estimate the mean using \eqref{eq:estmean1}, whereas the rule that is constructed using $q(x)$ requires integration of $f(x) = \rho(x) \cdot x$. All rules are constructed using $10^5$ samples and again the experiment is repeated 50 times to demonstrate general behavior of the rules. The results are gathered in Figure~\ref{subfig:interpolationerror}.

The rule constructed using samples from $\rho(x)$ immediately saturates to a level where the sampling error dominates, as no interpolation is used and the quadrature rule error vanishes since the error estimate from \eqref{eq:estmean1} is based on the mean (i.e.\ the first moment). The rule constructed using the prior $q(x)$ converges geometrically, as the PDF of the $\mathrm{Beta}(40, 60)$-distribution is a polynomial. A significant improvement is noticed by using our proposed adaptive proposal distribution $\rho_N(x)$. Initially a larger rate of convergence is obtained (up to approximately $N=10$), as the interpolant provides rapid convergence. For large $N$, the quadrature error starts dominating and the rate of convergence equals the rate of the rule constructed using $q(x)$.

This test case demonstrates that using samples that have been obtained directly from $\rho(x)$ is preferable, but this is intractable in practice. The adaptive proposal distribution has a clear advantage over using a rule that only incorporates the prior.

\subsubsection{Two-dimensional Genz test functions}
\label{subsubsec:peakypeaky}
\begin{figure}
	\begin{minipage}{\textwidth}
		\centering
		\small
		\includepgf{.6\textwidth}{.5\textwidth}{peaky-convergence.tikz}

		\caption{Convergence of the mean error determined with the new adaptive implicit quadrature rule or a quadrature rule determined using the prior. Equal colors refer to the same quadrature rule. Equal symbols refer to the same function.}
		\label{fig:peaky-convergence}
	\end{minipage}
	\vspace\intextsep
	\vspace\baselineskip

	\begin{minipage}{.3\textwidth}
		\centering
		\small
		\includepgf{\textwidth}{\textwidth}{peaky-1.tikz}
		\subcaption{$u_1$}
	\end{minipage}
	\hfill
	\begin{minipage}{.3\textwidth}
		\centering
		\small
		\includepgf{\textwidth}{\textwidth}{peaky-2.tikz}
		\subcaption{$u_2$}
	\end{minipage}
	\hfill
	\begin{minipage}{.3\textwidth}
		\centering
		\small
		\includepgf{\textwidth}{\textwidth}{peaky-3.tikz}
		\subcaption{$u_3$}
	\end{minipage}

	\caption{Examples of quadrature rules with $D = 65$ (exact for all polynomials up to degree 10) obtained by the adaptive implicit quadrature rule. The colors indicate the weights of the quadrature rule nodes. The nodes of a quadrature rule of the same polynomial degree, but constructed with respect to the prior, are depicted in gray.}
	\label{fig:peaky-quadrature-rules}
\end{figure}
The previously considered test case of Section~\ref{subsubsec:domination} demonstrates the advantages of using an adaptive quadrature rule. To illustrate the nodal placement of this rule and to further compare this rule to conventional quadrature rules constructed with respect to the prior, the calibration of three two-dimensional Genz test functions~\cite{Genz1984} with artificially generated data is considered. Therefore let $d = 2$, $\x^* = \trans{(1/2, 1/2)}$, and consider the following functions:
\begin{align}
	u_1(\x) &= \prod_{i=1}^d \left[ \frac{1}{4} + \left(x_i - \frac{1}{2}\right)^2 \right]^{-1}, &&\text{(Product Peak)} \\
	u_2(\x) &= \exp\left( -\sum_{i=1}^d \left|x_i - \frac{1}{2}\right| \right), &&\text{(C0 function)} \\
	u_3(\x) &= \begin{cases}
		0 &\text{if $x_1 > 3/5$ and $x_2 > 3/5$}, \\
		\exp\left( \sum_{i=1}^d x_i \right) &\text{otherwise}.
	\end{cases} &&\text{(Discontinuous function)}
\end{align}
These functions are specifically crafted such that they have difficult characteristics for numerical integration routines at or around $\x = \x^*$. The statistical model under consideration is
\begin{equation}
	\label{eq:statmodelgenz}
	z_k = u_g(\x) + \vareps, \text{ with } \vareps \sim \mathcal{N}(0, \sigma^2) \text{ and } \sigma^2 = 1/5, \text{ for $k = 1, \dots, 20$}.
\end{equation}
Here, $u_g$ (with $g = 1, 2, 3$) denotes one of the functions under consideration and $\mathbf{z} = (z_1, \dots, z_{20})^\textrm{T}$ are 20 ``measurements'' obtained by $z_k = u_g(\x^*) + n_k$ with $n_k$ a sample from $\mathcal{N}(0, \sigma^2)$. Hence the data is centered around the defining (``difficult'') characteristics of the test functions.

The goal is to infer $\x$ from \eqref{eq:statmodelgenz}, or equivalently, characterize the posterior $q(\x \mid \z)$ with a uniform prior and Gaussian likelihood defined in $\Omega = [0, 1]^d$, i.e.\ we use $\rho(\x) = q(\z \mid \x) \, q(\x)$ with
\begin{equation}
\label{eq:bayesgenz}
\begin{aligned}
	q(\x) &= \begin{cases}
		1 &\text{if $\x \in \Omega$}, \\
		0 &\text{otherwise},
	\end{cases} \\
	q(\z \mid \x) &\propto \exp\left[ -\frac{1}{2} \frac{\sum_{k=1}^{20} (u(\x) - z_k)^2}{\sigma^2} \right].
\end{aligned}
\end{equation}

The quantity used to measure the error of the quadrature rules is the maximum error of the mean of the posterior, i.e.\
\begin{equation}
	e_N = \begin{cases}
		\max\left(| \A_N \varphi_1 - \I^\rho \varphi_1 |, ~| \A_N \varphi_2 - \I^\rho \varphi_2 |\right) &\text{if $\A_N$ w.r.t.\ $\rho_N(\x)$}, \\
		\max\left(\left| \A_N \Bigl(\varphi_1 q(\z \mid x_1)\Bigr) - \I^\rho \varphi_1 \right|, ~\left| \A_N \Bigl(\varphi_2 q(\z \mid x_2)\Bigr) - \I^\rho \varphi_2 \right|\right) &\text{if $\A_N$ w.r.t.\ $q(\x)$},
	\end{cases}
\end{equation}
where $\A_N$ denotes a quadrature rule operator constructed using the adaptive method proposed in this article (i.e.\ with respect to $\rho_N$) or solely using the prior (i.e.\ with respect to $q(\x)$). In both cases the weights are scaled such that $\sum_k w_k = 1$. All quadrature rules are constructed using a linearly growing $D_i$, i.e.\ in iteration $i$ it is enforced that the quadrature rule integrates all $\varphi \in \Phi_i$ exactly (with $\Phi_i$ as introduced in Section~\ref{subsubsec:degreepos}). Recall that if $\x = (x_1, x_2)$ we have that $\varphi_1(\x) = x_1$ and $\varphi_2(\x) = x_2$.

There is randomness both in the generation of the quadrature rules and in the generation of the measurement data. Therefore, the error $e_N$ is determined 25 times and the obtained errors are averaged, which is denoted by $\bar{e}_N$. The convergence results are summarized in Figure~\ref{fig:peaky-convergence}. Examples of quadrature rules obtained by applying the proposed procedure are depicted in Figure~\ref{fig:peaky-quadrature-rules}. There are still oscillations visible, regardless of the averaging, though the error decay can clearly be observed.

Notice that the adaptive implicit quadrature rule outperforms the conventional quadrature rule in all three cases. The scatter plots of the nodes illustrate that the quadrature rules are indeed tailored specifically to the posterior under consideration. Since all rules have positive weights, estimates determined using these rules are unconditionally stable. This is an aforementioned advantage of the proposed method: the quadrature rules exhibit both exploitation (due to the interpolant) and exploration (due to the spread of the nodes) without suffering from instability.

\subsubsection{Five-dimensional Genz test functions}
\label{subsubsec:genz}
In this section, the calibration of the five-dimensional Genz test functions is considered. The setting is similar to the one of the previous section, i.e.\ let $d = 5$, $\x^* \in \mathbb{R}^d$ with $x^*_k = 1/2$, and consider the following functions:
\begin{align}
	u_1(\x) &= \cos\left(2\pi b_1 + \sum_{i=1}^d a_i x_i\right), &&\text{(Oscillatory)} \\
	u_2(\x) &= \prod_{i=1}^d \left(a_i^{-2} + (x_i - b_i)^2\right)^{-1}, &&\text{(Product Peak)} \\
	u_3(\x) &= \left(1 + \sum_{i=1}^d a_i x_i\right)^{-(d+1)}, &&\text{(Corner Peak)} \\
	u_4(\x) &= \exp\left(- \sum_{i=1}^d a_i^2 (x_i - b_i)^2 \right), &&\text{(Gaussian)} \\
	u_5(\x) &= \exp\left(- \sum_{i=1}^d a_i |x_i - b_i|\right), &&\text{($C_0$ function)} \\
	u_6(\x) &= \begin{cases}
		0 &\text{if $x_1 > b_1$ or $x_2 > b_2$}, \\
		\exp\left(\sum_{i=1}^d a_i x_i\right) &\text{otherwise}.
	\end{cases} &&\text{(Discontinuous)}
\end{align}
The vectors $\mathbf{a} = \trans{(a_1, \dots, a_d)}$ and $\mathbf{b} = \trans{(b_1, \dots, b_d)}$ are shape and translation vectors, respectively, and are chosen randomly in this test case. The elements of $\mathbf{a}$ enlarge the defining ``difficult'' property of the function and $\mathbf{b}$ translates the function in space. The statistical model under consideration is the same model considered in the previous section, see \eqref{eq:statmodelgenz} and \eqref{eq:bayesgenz}. Notice that, due to the random nature of $\mathbf{b}$, the probability that the difficult property of the functions is around $\mathbf{x} = \mathbf{x}^*$ is small.

\begin{figure}
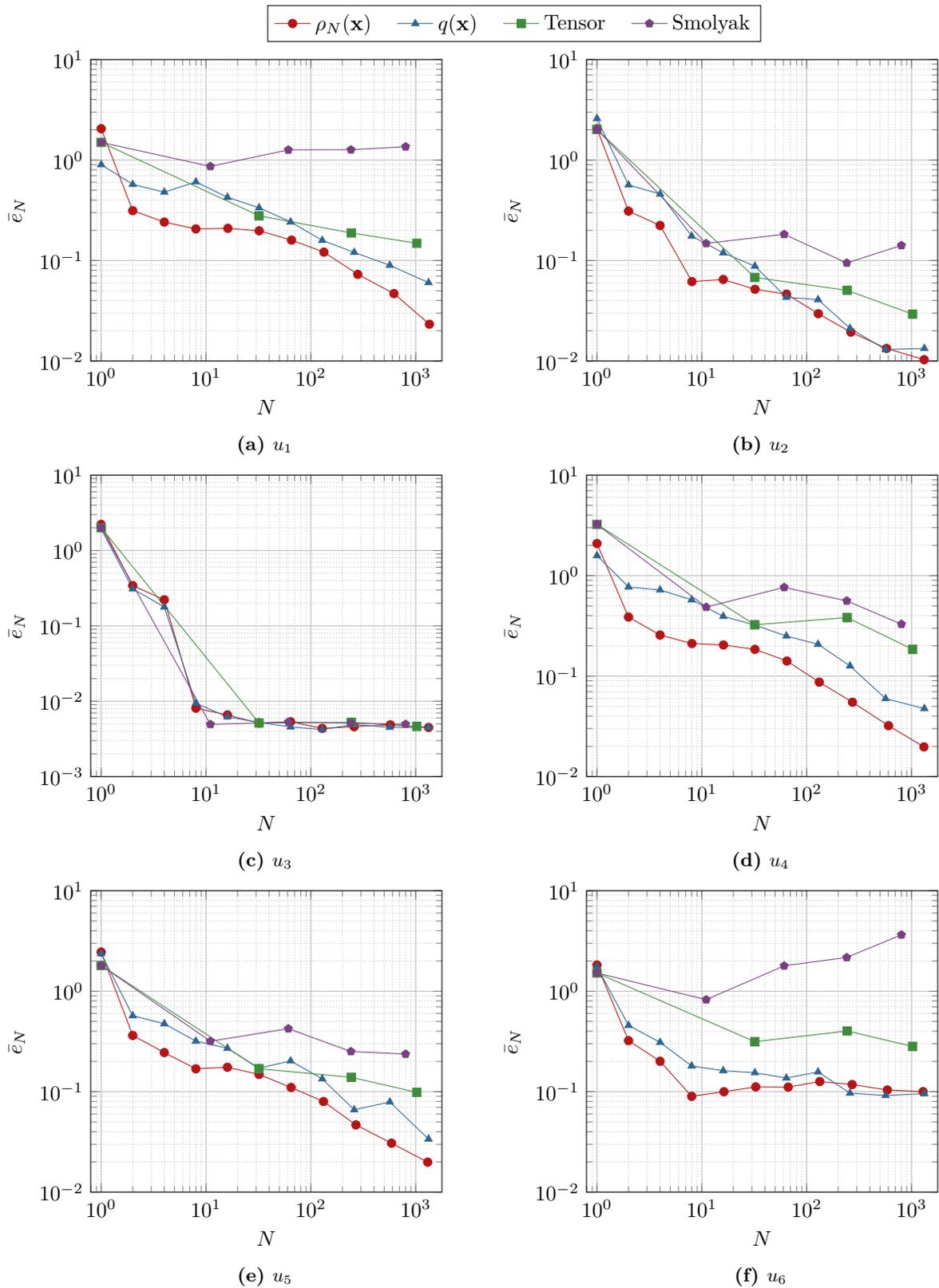

	\centering
	\pgfplotslegendfromname{genz-legend}

	\begin{minipage}[t]{.5\textwidth}
		\centering
		\includepgf{.9\textwidth}{.8\textwidth}{genz-1.tikz}
		\subcaption{$u_1$}
	\end{minipage}%
	\begin{minipage}[t]{.5\textwidth}
		\centering
		\includepgf{.9\textwidth}{.8\textwidth}{genz-2.tikz}
		\subcaption{$u_2$}
	\end{minipage}

	\begin{minipage}[t]{.5\textwidth}
		\centering
		\includepgf{.9\textwidth}{.8\textwidth}{genz-3.tikz}
		\subcaption{$u_3$}
	\end{minipage}%
	\begin{minipage}[t]{.5\textwidth}
		\centering
		\includepgf{.9\textwidth}{.8\textwidth}{genz-4.tikz}
		\subcaption{$u_4$}
	\end{minipage}

	\begin{minipage}[t]{.5\textwidth}
		\centering
		\includepgf{.9\textwidth}{.8\textwidth}{genz-5.tikz}
		\subcaption{$u_5$}
	\end{minipage}%
	\begin{minipage}[t]{.5\textwidth}
		\centering
		\includepgf{.9\textwidth}{.8\textwidth}{genz-6.tikz}
		\subcaption{$u_6$}
	\end{minipage}

	\caption{Convergence of the mean of the posterior determined using three different quadrature rules. The six functions under consideration are the Genz test functions.}
	\label{fig:genz}
\end{figure}

Four quadrature rules are considered to determine the mean of the posterior. Firstly, the proposed quadrature rule with a nearest neighbor interpolant is used. At each iteration the exactness of the quadrature rule is doubled, starting with $D_0 = 2^0$ up to $D_{10} = 2^{10}$. In other words, subsequent quadrature rules integrate larger numbers of polynomials (recall that $\Phi_D$ denotes the space of $D+1$ polynomials, sorted graded lexicographically). Secondly, an implicit quadrature rule without interpolation is used, generated using solely the prior. Finally, to compare the results with conventional quadrature rules, a tensor grid and a Smolyak sparse grid generated with a Clenshaw--Curtis quadrature rule are used. A linearly growing sequence of Clenshaw--Curtis quadrature rules is considered, so the tensor grid and the Smolyak sparse grid do not form a sequence of nested quadrature rules. This is done to keep the number of nodes of the multivariate quadrature rules tractable. We do not present the results obtained with Markov chain Monte Carlo, since the number of model evaluations considered here is too small to obtain a ``burned-in'' sequence.

For each quadrature rule, the experiment is repeated 25 times and the reported errors are averaged. For each experiment, the vectors $\mathbf{a}$ and $\mathbf{b}$ are drawn randomly from the five-dimensional unit hypercube and $\mathbf{a}$ is subsequently scaled such that $\| \mathbf{a} \|_2 = 5/2$. Moreover, for each experiment different samples are used to generate the implicit quadrature rules. To introduce randomness in the tensor grid and Smolyak sparse grid, the ``exact'' value $\x^*$ is randomized in the domain. This \emph{improves} the results of calibration with the Clenshaw--Curtis rules, as it otherwise would have a node at $\x = \x^*$ and therefore significantly overestimate the moments. The quantity used to measure the error of the quadrature rules is the same as in the previous test case, i.e.\ $e_N$ is the maximum error of all first-order raw moments (which are $\varphi_1, \dots, \varphi_5$). The ``exact'' mean is determined using a Monte Carlo sample from the prior. The number of samples used is the same number used to generate the implicit quadrature rules, such that the sampling error can be observed well (in that case all errors saturate). The obtained results are depicted in Figure~\ref{fig:genz}.

The first four test functions are analytic, so for these test functions the quadrature rules perform best. The flattening of the error of the third test function occurs due to the sampling error of the \emph{exact} mean (and not due to the sampling error of the quadrature rule). Moreover, the results of the third test function are best, since the corner peak varies only near a corner (hence the name) and therefore yields an almost uniform posterior.

The fifth Genz test function is only continuous (and not smooth), so should theoretically yield inferior performance. However, it is not differentiable only at a single point in the five-dimensional domain, so it is very unlikely that that point is of importance for the accuracy of the quadrature rules.

The sixth Genz test function is discontinuous and all rules do not yield a rapidly converging estimate. This is expected, since the quadrature rules are based on approximation using smooth functions. Notice that the convergence of the proposed adaptive procedure is significantly worse than in the two-dimensional case considered in Section~\ref{subsubsec:peakypeaky}, since detecting the discontinuity in the posterior using nearest neighbor interpolation in five dimensions requires significantly more model evaluations than in two dimensions. However, the accuracy remains similar to the other quadrature rules.

Overall, the adaptive quadrature rule proposed here consistently performs better or on the same level as the other quadrature rules. This demonstrates that the proposed procedure works for non-informative posteriors, since it is rarely the case that the posterior of the Genz test functions in conjunction with the randomized parameters is informative. Intuitively one would expect that piecewise interpolation is in these cases not necessarily beneficial, but the results demonstrate that it does not deteriorate the performance.

The Smolyak sparse grid performs significantly worse than the other quadrature rules for $u_1$ and $u_6$. This happens due to the negative weights of the sparse grid, or more specifically, due to the large positive weights present in the Smolyak sparse grid. These large weights significantly deteriorate the accuracy for a small number of the 25 runs, which is reflected in the averaged errors presented in Figure~\ref{fig:genz}. Notice that this effect can be significantly alleviated with an adaptive Smolyak sparse grid, as considered by \citet{Schillings2013}, though negative weights will remain present.

\subsubsection{d-dimensional Genz test functions}
\label{subsubsec:genzdim}
\begin{figure}
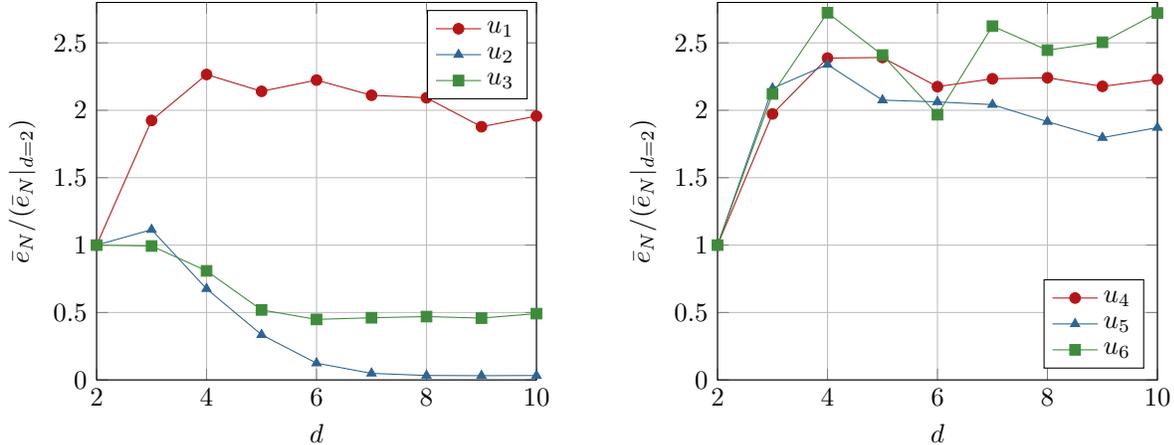

	\begin{minipage}{.5\textwidth}
		\centering
		\includepgf{.9\textwidth}{.8\textwidth}{genz-dim-123.tikz}
	\end{minipage}%
	\begin{minipage}{.5\textwidth}
		\centering
		\includepgf{.9\textwidth}{.8\textwidth}{genz-dim-456.tikz}
	\end{minipage}

	\caption{The integration error of the mean, averaged over 100 runs, considering the Genz test functions for varying dimension.}
	\label{fig:genz-dim}
\end{figure}
It is well-known that the convergence rate of quadrature rules based on polynomial approximation deteriorates in larger dimensional spaces. This follows among others from the Lebesgue inequality, see \eqref{eq:lebesgue}, since approximating multivariate functions requires larger numbers of basis polynomials. It is not straightforward to precisely assess the effect of the dimension on the accuracy of the quadrature rule, since the properties of the statistical model and the Genz test functions also depend on the dimension. However, if the accuracy of the quadrature rule deteriorates rapidly, the quadrature rule error will dominate the overall integration error, which is the quantity reported in this work.

Therefore, the Genz test functions are reconsidered in combination with the previously discussed statistical model, see \eqref{eq:bayesgenz}. The key difference is that the experiment is repeated for different dimensions $d$. Again the mean integration error $\bar{e}_N$ is computed, but contrary to the previous cases the error is determined by averaging over 100 experiments, since the error varies significantly less (and otherwise noise would pollute the results). The mean integration error is computed using a quadrature rule of degree 20, which is determined iteratively using a linearly growing sequence $D_i = i$ (hence 20 iterations of the algorithm are necessary to construct a quadrature rule). All other parameters are chosen as in Section~\ref{subsubsec:genz}. The results are reported in Figure~\ref{fig:genz-dim}, where the error is scaled with the integration error of the bivariate case (i.e.\ $d=2$).

The effect of the dimension on the statistical model and the Genz test functions is clearly visible. For $u_2$ and $u_3$, the error is \emph{decreasing} for increasing dimension, which conflicts with the theoretical expectations. This is due to the region that encompasses the defining properties of the product peak (i.e.\ $u_2$) and corner peak (i.e.\ $u_3$), which becomes smaller in multivariate spaces. A similar effect is also observed in the results of Section~\ref{subsubsec:genz}. The other functions yield a rapidly increasing error for $d \leq 5$, but the error stagnates (or even decreases slightly) for larger $d$. Again, this is the result of the varying properties of the test functions and the statistical model.

Even though it is not straightforward to isolate the error of the quadrature rule, this test case demonstrates that the growth of the quadrature error for increasing dimension does not dominate the growth of the integration error. In other words, the effect of dimension on the quadrature rule is significantly less than the effect of dimension on the test functions and the statistical model. Therefore the quadrature rules yield a viable approach for Bayesian prediction in multivariate spaces.

\subsection{Transonic flow over an airfoil}
\label{subsec:rae2822}
The applicability of the approach to complex test cases is demonstrated by considering the transonic flow over the RAE2822 airfoil. The problem under consideration is the calibration of the closure parameters of the Spalart--Allmaras turbulence model using wind tunnel measurements. It is computationally expensive to determine an accurate numerical solution of this problem, so the usage of efficient calibration approaches is of importance.

This section is split into three parts. Firstly, the problem of modeling the transonic flow over an airfoil and the measurement data under consideration is briefly introduced in Section~\ref{subsubsec:rans}. Secondly, the statistical model under consideration that relates the measurement data to the model output is discussed in Section~\ref{subsubsec:stats}. The obtained Bayesian predictions are presented in Section~\ref{subsubsec:results}.

\subsubsection{Modeling the flow over an airfoil}
\label{subsubsec:rans}
The flow over the airfoil under consideration is modeled using the Reynolds-averaged Navier--Stokes (RANS) equations, which only model the time-averaged mean flow over the airfoil. These equations do not form a closed system and require a closure model. For this purpose, the Boussinesq hypothesis is used, which introduces an eddy viscosity, which is subsequently modeled using a turbulence model. The details of this derivation are out of the scope of this article and we refer the interested reader to~\citet{Wilcox1998}. The turbulence model under consideration is the Spalart--Allmaras model~\cite{Spalart1992}, which is commonly used to model the flow over an airfoil. Various variants of the model exists, though in this work the original most commonly used variant from~\citet{Spalart1992} is considered. The model has 10 constants, which are typically defined as follows:
\begin{center}
	\begin{tabular}{r l | r l | r l}
		\hline
		$\sigma$ & $2/3$ & $C_{b1}$ & $0.1355$ & $C_{b2}$ & $0.622$ \\
		$\kappa$ & $0.41$ & $C_{w2}$ & $0.3$ & $C_{w3}$ & $2.0$ \\
		$C_{t3}$ & $1.2$ & $C_{t4}$ & $0.5$ & $C_{v1}$ & $7.1$  \\
		\hline
	\end{tabular}
\end{center}
Here $C_{w1}$ is omitted, since it depends on the other coefficients by the following expression:
\begin{equation}
	C_{w1} = C_{b1} / \kappa^2 + (1+C_{b2}) / \sigma.
\end{equation}

The closure coefficients of turbulence models are usually obtained by fitting the model using canonical test cases, for which analytical or accurate numerical solutions are available. There is no clear physical interpretation of many of these coefficients. A more structured approach to obtain values for these coefficients is to calibrate these coefficients statistically using Bayesian model calibration~\cite{Edeling2014-2,Cheung2011,Bos2018}, which has been discussed in this article. The posterior can be used to infer credible intervals of the parameters and can be propagated to infer Bayesian predictions.

A major advantage of the approach discussed in this article is that the determined quadrature rule can directly be used to infer prediction of the flow incorporating the uncertainty of the closure coefficients. The costly model evaluations have already been performed during the construction of the rule. No Markov chain Monte Carlo routines are necessary to obtain the predictions.

In this article, SU2~\cite{Palacios2014} is employed to numerically solve the RANS equations. It uses a second-order finite volume discretization and has support for the Spalart--Allmaras turbulence model. Moreover, it has been tested extensively to model the flow over an airfoil and the turbulence model parameters can be made configurable due to its freely available source code. As mentioned before, the airfoil under consideration is the RAE2822. Measurements of this airfoil for various values of the Reynolds number, Mach number, and angle of attack are freely available~\cite{Cook1979}. We consider Case 6 (see~\citet{Cook1979} for the other cases), with Reynolds number $6.5 \cdot 10^6$, Mach number $0.725$, and angle of attack $2.92^\circ$. The measurements encompass $n_z = 103$ point measurements of the pressure coefficient on the surface of the airfoil. To use these values in conjunction with SU2, we correct these values to incorporate wind tunnel effects~\cite{Slater2000}. The flow over the RAE2822 airfoil is transonic in this case and there is shock formation with a supersonic regime on the upper part of the airfoil.

\subsubsection{Statistical model}
\label{subsubsec:stats}
The quantity of interest used for the calibration is the pressure coefficient on the surface of the airfoil. Following~\citet{Edeling2014-2}, the parameters considered for calibration are $\theb = [\kappa,\allowbreak \sigma,\allowbreak C_{b1},\allowbreak C_{b2},\allowbreak C_{v1},\allowbreak C_{w2},\allowbreak C_{w3}]^\textrm{T}$. Let $s \in [0, 2]$ be the spatial parameter that runs from the trailing edge in anticlockwise direction over the airfoil. Then $u(\theb; s)$ denotes the mapping that yields the pressure coefficient at $s$ using the turbulence parameters $\theb$. We denote the calibration parameters by $\theb$ (instead of $\x$) to avoid confusion between $\theb$ and the usual spatial coordinates $x$ and $y$. A single evaluation of SU2 yields the pressure coefficient at all spatial coordinates in the discrete mesh for given $\theb$ and we use linear interpolation to obtain the pressure coefficient at the measurement locations.

The statistical model under consideration incorporates both model uncertainty, denoted by $\delta$, and measurement error, denoted by $\vareps$. It is, following~\citet{Kennedy2001}, as follows:
\begin{equation}
	\label{eq:statmodelrae2822}
	 z_k = u(\theb; s_k) + \delta(s_k) + \vareps_k, \text{ with } \vareps_k \sim \mathcal{N}(0, \sigma^2) \text{ for } k = 1, \dots, n_z.
\end{equation}
Here, $z_k$ denotes the measured pressure coefficient at location $s_k$. The random process $\delta(s) \sim \mathcal{N}(0, \Cov(s, s' \mid A, l))$ models the discrepancy between the model and the truth by means of a Gaussian process with zero mean and squared exponential covariance:
\begin{equation}
	\Cov(s, s' \mid A, l) = A \exp\left[ - \left( \frac{s - s'}{L 10^l} \right)^2 \right].
\end{equation}
The values of $A$ and $l$ are hyperparameters and are being inferred from calibration, whereas $L$ is a fixed parameter whose value is provided later. Measurement error is modeled by $\vareps_k$, which is Gaussian distributed with known standard deviation $\sigma = 0.01$. This value is based on the measurement errors reported by~\citet{Cook1979}.

The model from \eqref{eq:statmodelrae2822} yields the following likelihood: 
\begin{equation}
	\label{eq:likelihoodrae2822}
	q(\z \mid \theb, A, l) \propto \exp \left[ -\frac{1}{2} \trans{\mathbf{d}} K^{-1} \mathbf{d} \right],
\end{equation}
with $\mathbf{d}$ the misfit and $K$ the covariance matrix, i.e.
\begin{align}
	d_k &= z_k - u(\theb; s_k), \\
	K &= \Sigma + C, \text{ with } C_{i,j} = \Cov(s_i, s_j \mid A, l) \text{ and } \Sigma = \sigma I.
\end{align}

The prior of the turbulence closure parameters is as follows:
\begin{center}
	\begin{tabular}{r l}
		\hline
		$\kappa$ & $[0.205, 0.615]$ \\
		$\sigma$ & $[1/3, 1]$ \\
		$C_{b1}$ & $[0.0678, 0.2033]$ \\
		$C_{b2}$ & $[0.311, 0.933]$ \\
		$C_{v1}$ & $[3.55, 10.65]$ \\
		$C_{w2}$ & $[0.2, 0.4]$ \\
		$C_{w3}$ & $[1, 3]$ \\
		\hline
	\end{tabular}
\end{center}
These values are such that they encapsulate all commonly used variants of the turbulence model and are such that SU2 converges for all possible combinations.

The hyperparameters $A$ and $l$ are also being calibrated and therefore also require a prior. These are respectively $A \sim \mathcal{U}(0, 0.01)$ and $l \sim \mathcal{U}(0, 1)$. We choose $L = 0.01$, such that the maximal covariance length is significantly larger than the length of a single numerical cell on the surface of the airfoil. The computational model does not depend on $A$, $l$, and $L$, so given $\theb$ the distribution of these parameters is fully known analytically.

By combining the prior with the likelihood from \eqref{eq:likelihoodrae2822}, Bayes' law yields the posterior:
\begin{equation}
	q(\theb, A, l \mid \mathbf{z}) \propto q(\mathbf{z} \mid \theb, A, l) \, q(\theb) \, q(A) \, q(l).
\end{equation}
We are not interested in the exact distribution of the hyperparameters and the evidence can be neglected by rescaling the quadrature rule weights. Hence we apply the adaptive implicit quadrature rule to the following distribution:
\begin{equation}
	\rho(\theb) = \iint q(\mathbf{z} \mid \theb, A, l) \, q(\theb) \, q(A) \, q(l) \dd l \dd A.
\end{equation}
Afterwards, the quadrature rule weights are scaled such that $\sum_{k=0}^N w_k = 1$ and the rule is used to infer Bayesian predictions, e.g.\ as derived in \eqref{eq:prediction}:
\begin{equation}
	\mathbb{E}[\widehat{\z} \mid \z](s) = \int_\Omega u(\theb, s) \, q(\theb \mid \z) \dd \theb \approx \sum_{k=0}^N w_k \, u(\theb_k, s), \text{ with } q(\theb \mid \z) \propto \rho(\theb).
\end{equation}

\subsubsection{Results}
\label{subsubsec:results}
\begin{figure}
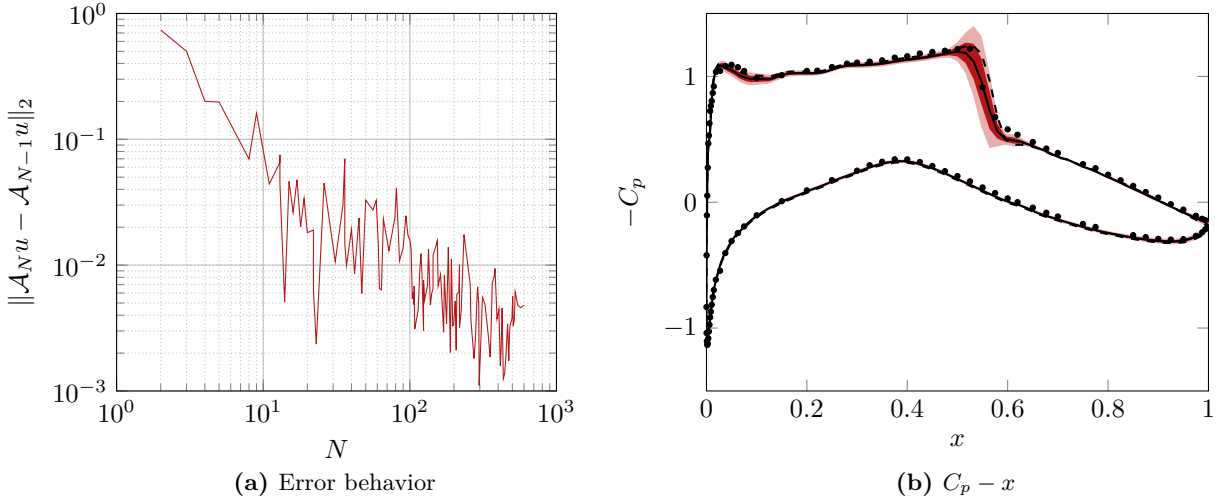

	\begin{minipage}[t]{.5\textwidth}
		\centering
		\includepgf{.9\textwidth}{.8\textwidth}{cp-convergence.tikz}
	\end{minipage}%
	\begin{minipage}[t]{.5\textwidth}
		\centering
		\includepgf{\textwidth}{.8\textwidth}{cp-x.tikz}
	\end{minipage}

	\begin{minipage}{.5\textwidth}
		\subcaption{Error behavior}
		\label{fig:cp-convergence}
	\end{minipage}%
	\begin{minipage}{.5\textwidth}
		\subcaption{$C_p - x$}
		\label{fig:cp-x}
	\end{minipage}

	\caption{\textit{Left:} convergence of consecutive differences, \textit{Right:} Mean and standard deviation of $C_p-x$ determined using the proposed approach (solid line and red area), compared to using canonical coefficients (dashed) and measurement data (scattered)}
	\label{fig:cp-x-convergence}
\end{figure}
The quadrature rule constructed in this test case has polynomial degree 3, which constitutes a polynomial space with 120 basis polynomials. This number of polynomials is a good balance between accuracy in the predictions and overall running time. The quadrature rules are constructed with linearly increasing exactness, so at the $i$-th iteration a quadrature rule of degree $i$ is constructed that incorporates all available model evaluations (up to iteration $i-1$).

The obtained rule incorporates an approximation of the posterior. There is no exact solution available, so it is not straightforward to assess the accuracy of the quadrature rule. For this purpose, each iteration the mean of the pressure coefficient is determined and these are consecutively compared using the 2-norm, i.e.\ if $N_i$ is the number of nodes of the quadrature rule after $i$ iterations, then the error measure is
\begin{equation}
	e_{N_i} = \| \A_{N_i} u - \A_{N_{i-1}} u \|_2.
\end{equation}
We will denote this with a little abuse of notation simply as
\begin{equation}
	e_N = \| \A_N u - \A_{N-1} u \|_2.
\end{equation}
The obtained errors are depicted in Figure~\ref{fig:cp-convergence}.

The quadrature rule clearly yields a converging mean and the error converges approximately linearly to zero. The oscillations are the result of the random sampling of the proposal distributions in combination with the fact that this test case was run only once. Obtaining a smoother decay would require repeating the experiment various times and averaging the obtained error (as done in the previous test cases), but this is too computationally costly in this case. In practical computations, the error decay can be assessed by regressing the convergence rate by means of least squares (i.e.\ by fitting analytic error decay).

Each iteration SU2 is evaluated for each additional quadrature rule node, so obtaining Bayesian predictions of the pressure coefficient is a straightforward post-processing step. The obtained predictions of the pressure coefficient are depicted in Figure~\ref{fig:cp-x}. The measurement data, uncertainty bounds, and the deterministic pressure coefficient determined using the canonical turbulence coefficients are also depicted. It is clearly visible that the largest uncertainty occurs near the shock on top of the airfoil, which is in line with existing results on uncertain transonic flows around the RAE2822 airfoil, either using uncorrelated prescribed distributions~\cite{Witteveen2009a,Pisaroni2017} or calibrated parameters in a Bayesian setting~\cite{Bos2018}. Notice that the number of evaluations of SU2 needed to infer these predictions is of similar order of magnitude as is common for uncertainty propagation with \emph{prescribed} distributions. It is a significant reduction compared to commonly used Markov chain Monte Carlo methods~\cite{Edeling2014-2,Edeling2014,Cheung2011}.

\begin{figure}
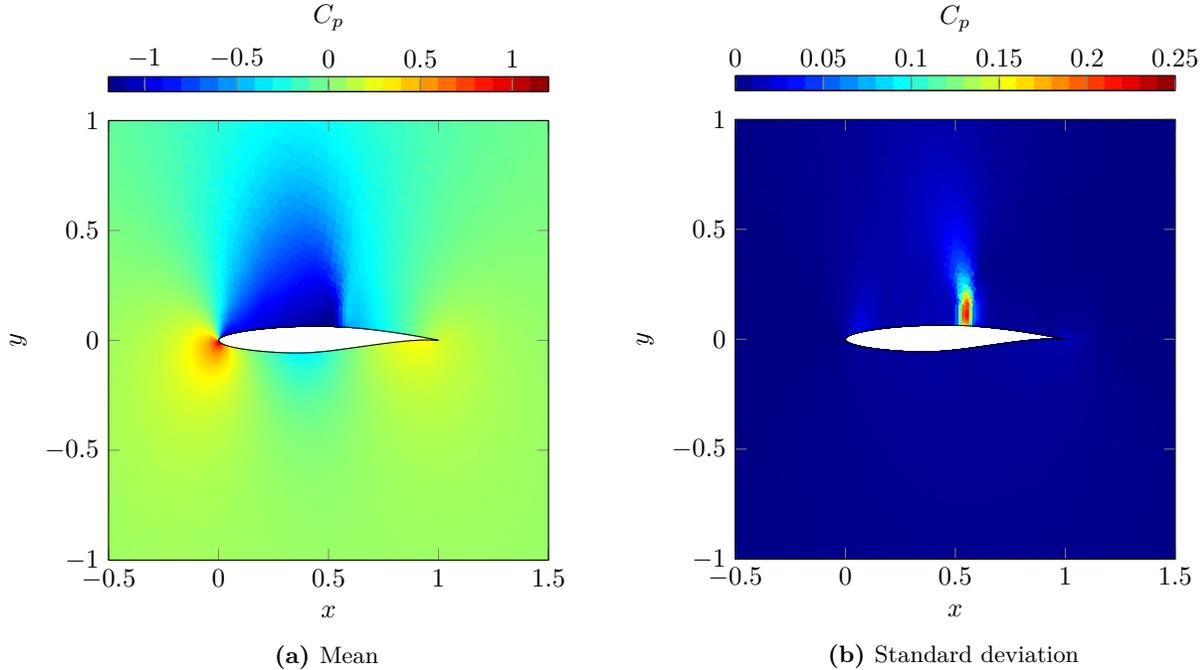

	\begin{minipage}{.5\textwidth}
		\centering
		\includepgf{.9\textwidth}{.9\textwidth}{cp-mean.tikz}
		\subcaption{Mean}
	\end{minipage}
	\begin{minipage}{.5\textwidth}
		\centering
		\includepgf{.9\textwidth}{.9\textwidth}{cp-stddev.tikz}
		\subcaption{Standard deviation}
	\end{minipage}

	\caption{Predictions of pressure coefficients in vicinity of the airfoil using calibrated coefficients.}
	\label{fig:cp-mean-stddev}
\end{figure}

SU2 yields the full flow field around the airfoil and therefore it is also possible to predict the statistical moments of the pressure coefficient away from the airfoil surface, as depicted in Figure~\ref{fig:cp-mean-stddev}. We want to emphasize the importance of positive weights, since these ensure that the prediction of the variance is positive. Notice that the pressure coefficients depicted in these figures (i.e.\ away from the airfoil) have \emph{not} been used in the calibration process, as only measurement data on the airfoil surface is available. So technically it not evident whether these results are correct or not. It is not straightforward to make general claims about this, since the specific properties (such as smoothness) of the model under consideration determine to a large extent whether it is viable to infer predictions of quantities using parameters that have been calibrated with a different quantity.

\section{Conclusions}
\label{sec:conclusion}
In this article a new methodology is proposed to construct quadrature rules with positive weights and high degree for the purpose of inferring Bayesian predictions. It consists of two steps. Firstly a quadrature rule is constructed using an approximation of the posterior. For this step the so-called implicit quadrature rule has been used. Secondly the approximation of the posterior is refined using all available model evaluations. For the second step nearest neighbor interpolation has been used. It has been demonstrated rigorously that our approach yields a quadrature rule whose error behaves like a quadrature rule with respect to the exact posterior for increasing number of nodes.

The performance of the quadrature rule has been demonstrated by calibrating the Genz test functions in a Bayesian framework. The results demonstrate that the quadrature rules consistently outperform (or perform as well as) conventional numerical integration approaches. The performance gains are most significant if the posterior is very informative, e.g.\ if its standard deviation is small. If, on the other hand, the posterior resembles approximately the prior, the performance is comparable to a quadrature rule determined using solely the prior.

The applicability of the approach to an expensive fluid dynamics model has been demonstrated by using the quadrature rule to calibrate the transonic flow over the RAE2822 airfoil. It has been demonstrated numerically that estimates of the rule converge and predictions of the mean and standard deviation of the pressure coefficient have been inferred. The results are in line with existing research, though significantly less model evaluations are necessary compared to existing Bayesian calibration approaches.

There are various opportunities to further extend the approach. The nearest neighbor interpolation procedure does not leverage smoothness in either the distribution or the posterior. A major opportunity is replacing nearest neighbor interpolation by a more efficient interpolation approach, which could allow for more rapid convergence. This is however far from trivial, since our method requires that the obtained interpolant is a distribution.

\section*{Acknowledgments}
This research is part of the Dutch EUROS program, which is supported by NWO domain Applied and Engineering Sciences and partly funded by the Dutch Ministry of Economic Affairs.

\small
\bibliographystyle{plainnatnourl}
\bibliography{literature.bib}

\begin{thebibliography}{67}
\providecommand{\natexlab}[1]{#1}
\providecommand{\url}[1]{\texttt{#1}}
\expandafter\ifx\csname urlstyle\endcsname\relax
  \providecommand{\doi}[1]{doi: #1}\else
  \providecommand{\doi}{doi: \begingroup \urlstyle{rm}\Url}\fi

\bibitem[Ahn et~al.(2012)Ahn, Korattikara, and Welling]{Ahn2012}
S.~Ahn, A.~Korattikara, and M.~Welling.
\newblock Bayesian posterior sampling via stochastic gradient {Fisher} scoring.
\newblock In \emph{Proceedings of the 29\textsuperscript{th} International
  Conference on Machine Learning}, 2012.

\bibitem[Biau and Devroye(2015)]{Biau2015}
G.~Biau and L.~Devroye.
\newblock \emph{Lectures on the Nearest Neighbor Method}.
\newblock Springer International Publishing, 2015.
\newblock \doi{10.1007/978-3-319-25388-6}.

\bibitem[Bilionis and Zabaras(2013)]{Bilionis2013}
I.~Bilionis and N.~Zabaras.
\newblock Solution of inverse problems with limited forward solver evaluations:
  a {Bayesian} perspective.
\newblock \emph{Inverse Problems}, 30\penalty0 (1):\penalty0 015004, 2013.
\newblock \doi{10.1088/0266-5611/30/1/015004}.

\bibitem[Birolleau et~al.(2014)Birolleau, Po\"ette, and Lucor]{Birolleau2014}
A.~Birolleau, G.~Po\"ette, and D.~Lucor.
\newblock Adaptive {B}ayesian inference for discontinuous inverse problems,
  application to hyperbolic conservation laws.
\newblock \emph{Communications in Computational Physics}, 16\penalty0
  (1):\penalty0 1--34, 2014.
\newblock \doi{10.4208/cicp.240113.071113a}.

\bibitem[van~den Bos et~al.(2017)van~den Bos, Koren, and Dwight]{Bos2016b}
L.~M.~M. van~den Bos, B.~Koren, and R.~P. Dwight.
\newblock Non-intrusive uncertainty quantification using reduced cubature
  rules.
\newblock \emph{Journal of Computational Physics}, 332:\penalty0 418--445,
  2017.
\newblock \doi{10.1016/j.jcp.2016.12.011}.

\bibitem[van~den Bos et~al.(2020{\natexlab{a}})van~den Bos, Sanderse,
  Bierbooms, and van Bussel]{Bos2018}
L.~M.~M. van~den Bos, B.~Sanderse, W.~A. A.~M. Bierbooms, and G.~J.~W. van
  Bussel.
\newblock Bayesian model calibration with interpolating polynomials based on
  adaptively weighted {Leja} nodes.
\newblock \emph{Communications in Computational Physics}, 27\penalty0
  (1):\penalty0 33--69, 2020{\natexlab{a}}.
\newblock \doi{10.4208/cicp.oa-2018-0218}.

\bibitem[van~den Bos et~al.(2020{\natexlab{b}})van~den Bos, Sanderse,
  Bierbooms, and van Bussel]{Bos2018b}
L.~M.~M. van~den Bos, B.~Sanderse, W.~A. A.~M. Bierbooms, and G.~J.~W. van
  Bussel.
\newblock Generating nested quadrature rules with positive weights based on
  arbitrary sample sets.
\newblock \emph{SIAM/ASA Journal on Uncertainty Quantification}, 8\penalty0
  (1):\penalty0 139--169, 2020{\natexlab{b}}.
\newblock \doi{10.1137/18M1213373}.

\bibitem[Botts et~al.(2011)Botts, H{\"{o}}rmann, and Leydold]{Botts2011}
C.~Botts, W.~H{\"{o}}rmann, and J.~Leydold.
\newblock Transformed density rejection with inflection points.
\newblock \emph{Statistics and Computing}, 23\penalty0 (2):\penalty0 251--260,
  2011.
\newblock \doi{10.1007/s11222-011-9306-4}.

\bibitem[Brandolini et~al.(2014)Brandolini, Choirat, and
  Colzani]{Brandolini2014}
L.~Brandolini, C.~Choirat, and L.~Colzani.
\newblock Quadrature rules and distribution of points on manifolds.
\newblock \emph{Annali della Scuola Normale Superiore - Classe di Scienze},
  2014.
\newblock ISSN 2036-2145.
\newblock \doi{10.2422/2036-2145.201103_007}.

\bibitem[Brass and Petras(2011)]{Brass2011}
H.~Brass and K.~Petras.
\newblock \emph{Quadrature Theory}, volume 178 of \emph{Mathematical Surveys
  and Monographs}.
\newblock American Mathematical Society, 2011.
\newblock \doi{10.1090/surv/178}.

\bibitem[Brutman(1996)]{Brutman1996}
L.~Brutman.
\newblock Lebesgue functions for polynomial interpolation-a survey.
\newblock \emph{Annals of Numerical Mathematics}, 4:\penalty0 111--128, 1996.

\bibitem[Buhmann(2000)]{Buhmann2000}
M.~D. Buhmann.
\newblock Radial basis functions.
\newblock \emph{Acta Numerica}, pages 1--38, 2000.
\newblock \doi{10.1007/978-3-540-69909-5_3}.

\bibitem[Butler et~al.(2017)Butler, Graham, Mattis, and Walsh]{Butler2017}
T.~Butler, L.~Graham, S.~Mattis, and S.~Walsh.
\newblock A measure-theoretic interpretation of sample based numerical
  integration with applications to inverse and prediction problems under
  uncertainty.
\newblock \emph{{SIAM} Journal on Scientific Computing}, 39\penalty0
  (5):\penalty0 A2072--A2098, 2017.
\newblock \doi{10.1137/16m1063289}.

\bibitem[Butler et~al.(2018)Butler, Jakeman, and Wildey]{Butler2018a}
T.~Butler, J.~D. Jakeman, and T.~Wildey.
\newblock Convergence of probability densities using approximate models for
  forward and inverse problems in uncertainty quantification.
\newblock \emph{{SIAM} Journal on Scientific Computing}, 40\penalty0
  (5):\penalty0 A3523--A3548, 2018.
\newblock \doi{10.1137/18m1181675}.

\bibitem[Caflisch(1998)]{Caflisch1998}
R.~E. Caflisch.
\newblock {Monte} {Carlo} and quasi-{Monte} {Carlo} methods.
\newblock \emph{{ANU}}, 7:\penalty0 1, 1998.
\newblock \doi{10.1017/s0962492900002804}.

\bibitem[Cheung et~al.(2011)Cheung, Oliver, Prudencio, Prudhomme, and
  Moser]{Cheung2011}
S.~H. Cheung, T.~A. Oliver, E.~E. Prudencio, S.~Prudhomme, and R.~D. Moser.
\newblock Bayesian uncertainty analysis with applications to turbulence
  modeling.
\newblock \emph{Reliability Engineering {\&} System Safety}, 96\penalty0
  (9):\penalty0 1137--1149, 2011.
\newblock \doi{10.1016/j.ress.2010.09.013}.

\bibitem[Cook et~al.(1979)Cook, McDonald, and Firmin]{Cook1979}
P.~H. Cook, M.~A. McDonald, and M.~C.~P. Firmin.
\newblock Aerofoil {RAE 2822} -- pressure distributions, and boundary layer and
  wake measurements.
\newblock In \emph{Experimental Data Base for Computer Program Assessment},
  number 138 in AGARD Advisory Report, chapter~6, pages A6--1 -- A6--77. North
  Atlantic Treaty Organization, 1979.

\bibitem[Derflinger et~al.(2010)Derflinger, H{\"{o}}rmann, and
  Leydold]{Derflinger2010}
G.~Derflinger, W.~H{\"{o}}rmann, and J.~Leydold.
\newblock Random variate generation by numerical inversion when only the
  density is known.
\newblock \emph{{ACM} Transactions on Modeling and Computer Simulation},
  20\penalty0 (4):\penalty0 1--25, 2010.
\newblock \doi{10.1145/1842722.1842723}.

\bibitem[Du et~al.(1999)Du, Faber, and Gunzburger]{Du1999}
Q.~Du, V.~Faber, and M.~Gunzburger.
\newblock Centroidal {Voronoi} tessellations: Applications and algorithms.
\newblock \emph{{SIAM} Review}, 41\penalty0 (4):\penalty0 637--676, 1999.
\newblock \doi{10.1137/s0036144599352836}.

\bibitem[Ebeida et~al.(2016)Ebeida, Mitchell, Swiler, Romero, and
  Rushdi]{Ebeida2016}
M.~S. Ebeida, S.~A. Mitchell, L.~P. Swiler, V.~J. Romero, and A.~A. Rushdi.
\newblock {POF}-darts: Geometric adaptive sampling for probability of failure.
\newblock \emph{Reliability Engineering {\&} System Safety}, 155:\penalty0
  64--77, 2016.
\newblock \doi{10.1016/j.ress.2016.05.001}.

\bibitem[Edeling et~al.(2014{\natexlab{a}})Edeling, Cinnella, and
  Dwight]{Edeling2014}
W.~N. Edeling, P.~Cinnella, and R.~P. Dwight.
\newblock Predictive {RANS} simulations via {B}ayesian model-scenario
  averaging.
\newblock \emph{Journal of Computational Physics}, 275:\penalty0 65--91,
  2014{\natexlab{a}}.
\newblock \doi{10.1016/j.jcp.2014.06.052}.

\bibitem[Edeling et~al.(2014{\natexlab{b}})Edeling, Cinnella, Dwight, and
  Bijl]{Edeling2014-2}
W.~N. Edeling, P.~Cinnella, R.~P. Dwight, and H.~Bijl.
\newblock {B}ayesian estimates of parameter variability in the $k-\epsilon$
  turbulence model.
\newblock \emph{Journal of Computational Physics}, 258:\penalty0 73--94,
  2014{\natexlab{b}}.
\newblock \doi{10.1016/j.jcp.2013.10.027}.

\bibitem[Edeling et~al.(2016)Edeling, Dwight, and Cinnella]{Edeling2016}
W.~N. Edeling, R.~P. Dwight, and P.~Cinnella.
\newblock Simplex-stochastic collocation method with improved scalability.
\newblock \emph{Journal of Computational Physics}, 310:\penalty0 301--328,
  2016.
\newblock \doi{10.1016/j.jcp.2015.12.034}.

\bibitem[Franke(1982)]{Franke1982}
R.~Franke.
\newblock Scattered data interpolation: tests of some methods.
\newblock \emph{Mathematics of Computation}, 38\penalty0 (157):\penalty0
  181--200, 1982.
\newblock \doi{10.1090/s0025-5718-1982-0637296-4}.

\bibitem[Fritsch and Carlson(1980)]{Fritsch1980}
F.~N. Fritsch and R.~E. Carlson.
\newblock Monotone piecewise cubic interpolation.
\newblock \emph{{SIAM} Journal on Numerical Analysis}, 17\penalty0
  (2):\penalty0 238--246, 1980.
\newblock \doi{10.1137/0717021}.

\bibitem[Gelman et~al.(2013)Gelman, Carlin, Stern, Dunson, Vehtari, and
  Rubin]{Gelman2013}
A.~Gelman, J.~B. Carlin, H.~S. Stern, D.~B. Dunson, A.~Vehtari, and D.~B.
  Rubin.
\newblock \emph{{B}ayesian Data Analysis}.
\newblock Chapman and Hall/CRC, third edition, 2013.

\bibitem[Genz(1984)]{Genz1984}
A.~Genz.
\newblock Testing multidimensional integration routines.
\newblock In \emph{Proceedings of International Conference on Tools, Methods
  and Languages for Scientific and Engineering Computation}, pages 81--94.
  Elsevier North--Holland, 1984.

\bibitem[Golub and Welsch(1969)]{Golub1969}
G.~H. Golub and J.~H. Welsch.
\newblock Calculation of {Gauss} quadrature rules.
\newblock \emph{Mathematics of Computation}, 23\penalty0 (106):\penalty0
  221--230, 1969.
\newblock \doi{10.1090/s0025-5718-69-99647-1}.

\bibitem[Grzelak et~al.(2018)Grzelak, Witteveen, Su{\'{a}}rez-Taboada, and
  Oosterlee]{Grzelak2018}
L.~A. Grzelak, J.~A.~S. Witteveen, M.~Su{\'{a}}rez-Taboada, and C.~W.
  Oosterlee.
\newblock The stochastic collocation {Monte} {Carlo} sampler: highly efficient
  sampling from `expensive' distributions.
\newblock \emph{Quantitative Finance}, pages 1--18, 2018.
\newblock \doi{10.1080/14697688.2018.1459807}.

\bibitem[Guessab and Schmeisser(2008)]{Guessab2008}
A.~Guessab and G.~Schmeisser.
\newblock Construction of positive definite cubature formulae and approximation
  of functions via {Voronoi} tessellations.
\newblock \emph{Advances in Computational Mathematics}, 32\penalty0
  (1):\penalty0 25--41, 2008.
\newblock \doi{10.1007/s10444-008-9080-9}.

\bibitem[Hastings(1970)]{Hastings1970}
W.~K. Hastings.
\newblock {Monte} {Carlo} sampling methods using {Markov} chains and their
  applications.
\newblock \emph{Biometrika}, 57\penalty0 (1):\penalty0 97--109, 1970.
\newblock \doi{10.1093/biomet/57.1.97}.

\bibitem[Hewitt and Hoeting(2019)]{Hewitt2019}
J.~Hewitt and J.~A. Hoeting.
\newblock Approximate {Bayesian} inference via sparse grid quadrature
  evaluation for hierarchical models.
\newblock \emph{ArXiv 1904.07270}, 2019.

\bibitem[Ibrahimoglu(2016)]{Ibrahimoglu2016}
B.~A. Ibrahimoglu.
\newblock Lebesgue functions and {Lebesgue} constants in polynomial
  interpolation.
\newblock \emph{Journal of Inequalities and Applications}, 2016\penalty0
  (1):\penalty0 93, 2016.
\newblock \doi{10.1186/s13660-016-1030-3}.

\bibitem[Jackson(1982)]{Jackson1982}
D.~Jackson.
\newblock \emph{Theory of Approximation (Colloquium Publications)}.
\newblock American Mathematical Society, 1982.
\newblock ISBN 978-0821810118.

\bibitem[Jakeman and Narayan(2018)]{Jakeman2017}
J.~D. Jakeman and A.~Narayan.
\newblock Generation and application of multivariate polynomial quadrature
  rules.
\newblock \emph{Computer Methods in Applied Mechanics and Engineering},
  338:\penalty0 134--161, 2018.
\newblock \doi{10.1016/j.cma.2018.04.009}.

\bibitem[Kennedy and O'Hagan(2001)]{Kennedy2001}
M.~C. Kennedy and A.~O'Hagan.
\newblock {B}ayesian calibration of computer models.
\newblock \emph{Journal of the Royal Statistical Society: Series B (Statistical
  Methodology)}, 63\penalty0 (3):\penalty0 425--464, 2001.
\newblock \doi{10.1111/1467-9868.00294}.

\bibitem[Le~Ma{\^{i}}tre and Knio(2010)]{LeMaitre2010}
O.~P. Le~Ma{\^{i}}tre and O.~M. Knio.
\newblock Spectral methods for uncertainty quantification.
\newblock \emph{Scientific Computation}, 2010.
\newblock ISSN 1434-8322.
\newblock \doi{10.1007/978-90-481-3520-2}.

\bibitem[Li and Marzouk(2014)]{Li2014}
J.~Li and Y.~M. Marzouk.
\newblock Adaptive construction of surrogates for the {Bayesian} solution of
  inverse problems.
\newblock \emph{SIAM Journal on Scientific Computing}, 36\penalty0
  (3):\penalty0 A1163--A1186, 2014.
\newblock ISSN 1095-7197.
\newblock \doi{10.1137/130938189}.

\bibitem[Luo et~al.(2009)Luo, Sun, and Wang]{Luo2009}
C.~Luo, J.~Sun, and Y.~Wang.
\newblock Integral estimation from point cloud in d-dimensional space.
\newblock In \emph{Proceedings of the 25\textsuperscript{th} annual symposium
  on Computational geometry - {SCG}~'09}. {ACM} Press, 2009.
\newblock \doi{10.1145/1542362.1542389}.

\bibitem[Marzouk and Xiu(2009)]{Marzouk2009}
Y.~Marzouk and D.~Xiu.
\newblock A stochastic collocation approach to {B}ayesian inference in inverse
  problems.
\newblock \emph{Communications in Computational Physics}, 6\penalty0
  (4):\penalty0 826--847, 2009.
\newblock \doi{10.4208/cicp.2009.v6.p826}.

\bibitem[Metropolis et~al.(1953)Metropolis, Rosenbluth, Rosenbluth, Teller, and
  Teller]{Metropolis1953}
N.~Metropolis, A.~W. Rosenbluth, M.~N. Rosenbluth, A.~H. Teller, and E.~Teller.
\newblock Equation of state calculations by fast computing machines.
\newblock \emph{The Journal of Chemical Physics}, 21\penalty0 (6):\penalty0
  1087, 1953.
\newblock \doi{10.1063/1.1699114}.

\bibitem[Narayan and Jakeman(2014)]{Narayan2014}
A.~Narayan and J.~D. Jakeman.
\newblock Adaptive {Leja} sparse grid constructions for stochastic collocation
  and high-dimensional approximation.
\newblock \emph{{SIAM} Journal on Scientific Computing}, 36\penalty0
  (6):\penalty0 A2952--A2983, 2014.
\newblock \doi{10.1137/140966368}.

\bibitem[Niederreiter(1992)]{Niederreiter1992}
H.~Niederreiter.
\newblock \emph{Random Number Generation and Quasi-{Monte} {Carlo} Methods}.
\newblock Society for Industrial \& Applied Mathematics ({SIAM}), 1992.
\newblock \doi{10.1137/1.9781611970081}.

\bibitem[Novak and Ritter(1999)]{Novak1999}
E.~Novak and K.~Ritter.
\newblock Simple cubature formulas with high polynomial exactness.
\newblock \emph{Constructive Approximation}, 15\penalty0 (4):\penalty0 499,
  1999.
\newblock ISSN 1432-0940.
\newblock \doi{10.1007/s003659900119}.

\bibitem[Palacios et~al.(2014)Palacios, Economon, Aranake, Copeland, Lonkar,
  Lukaczyk, Manosalvas, Naik, Padron, Tracey, Variyar, and
  Alonso]{Palacios2014}
F.~Palacios, T.~D. Economon, A.~Aranake, S.~R. Copeland, A.~K. Lonkar, T.~W.
  Lukaczyk, D.~E. Manosalvas, K.~R. Naik, S.~Padron, B.~Tracey, A.~Variyar, and
  J.~J. Alonso.
\newblock Stanford university unstructured ({SU}2): Analysis and design
  technology for turbulent flows.
\newblock In \emph{52\textsuperscript{nd} Aerospace Sciences Meeting}. American
  Institute of Aeronautics and Astronautics, 2014.
\newblock \doi{10.2514/6.2014-0243}.

\bibitem[Peherstorfer and Marzouk(2019)]{Peherstorfer2019}
B.~Peherstorfer and Y.~Marzouk.
\newblock A transport-based multifidelity preconditioner for {Markov} chain
  {Monte} {Carlo}.
\newblock \emph{Advances in Computational Mathematics}, pages 1--28, 2019.
\newblock \doi{10.1007/s10444-019-09711-y}.

\bibitem[Pisaroni et~al.(2017)Pisaroni, Nobile, and Leyland]{Pisaroni2017}
M.~Pisaroni, F.~Nobile, and P.~Leyland.
\newblock Continuation multi-level {Monte-Carlo} method for uncertainty
  quantification in turbulent compressible aerodynamics problems modeled by
  {RANS}.
\newblock Technical Report 10.2017, {{\'{E}}}cole Polytechnique
  F{\'{e}}d{\'{e}}rale de Lausanne, 2017.

\bibitem[Rasmussen and Williams(2006)]{Rasmussen2006}
C.~E. Rasmussen and C.~K.~I. Williams.
\newblock \emph{Gaussian Processes for Machine Learning (Adaptive Computation
  And Machine Learning)}.
\newblock The MIT Press, 2006.
\newblock ISBN 0-262-18253-X.

\bibitem[Rushdi et~al.(2017)Rushdi, Swiler, Phipps, D'Elia, and
  Ebeida]{Rushdi2017}
A.~A. Rushdi, L.~P. Swiler, E.~T. Phipps, M.~D'Elia, and M.~S. Ebeida.
\newblock {VPS}: Voronoi piecewise surrogate models for high-dimensional data
  fitting.
\newblock \emph{International Journal for Uncertainty Quantification},
  7\penalty0 (1):\penalty0 1--21, 2017.
\newblock \doi{10.1615/int.j.uncertaintyquantification.2016018697}.

\bibitem[Schillings and Schwab(2013)]{Schillings2013}
C.~Schillings and C.~Schwab.
\newblock Sparse, adaptive {Smolyak} quadratures for {Bayesian} inverse
  problems.
\newblock \emph{Inverse Problems}, 29\penalty0 (6):\penalty0 065011, 2013.
\newblock \doi{10.1088/0266-5611/29/6/065011}.

\bibitem[Schillings et~al.(2019)Schillings, Sprungk, and Wacker]{Wacker2017}
C.~Schillings, B.~Sprungk, and P.~Wacker.
\newblock On the convergence of the {Laplace} approximation and
  noise-level-robustness of {Laplace}-based {Monte} {Carlo} methods for
  {Bayesian} inverse problems.
\newblock \emph{ArXiV 1901.03958}, 2019.

\bibitem[Schwab and Stuart(2012)]{Schwab2012}
C.~Schwab and A.~M. Stuart.
\newblock Sparse deterministic approximation of {Bayesian} inverse problems.
\newblock \emph{Inverse Problems}, 28\penalty0 (4):\penalty0 045003, 2012.
\newblock \doi{10.1088/0266-5611/28/4/045003}.

\bibitem[Slater et~al.(2000)Slater, Dudek, and Tatum]{Slater2000}
J.~W. Slater, J.~C. Dudek, and K.~E. Tatum.
\newblock The {NPARC} alliance verification and validation archive.
\newblock Technical Report 2000-209946, NASA, 2000.

\bibitem[Smolyak(1963)]{Smolyak1963}
S.~Smolyak.
\newblock Quadrature and interpolation formulas for tensor products of certain
  classes of functions.
\newblock \emph{Soviet Mathematics, Doklady}, 4:\penalty0 240--243, 1963.

\bibitem[Spalart and Allmaras(1992)]{Spalart1992}
P.~R. Spalart and S.~R. Allmaras.
\newblock A one-equation turbulence model for aerodynamic flows.
\newblock In \emph{30\textsuperscript{th} Aerospace Sciences Meeting \&
  Exhibit}, number 92-0439, pages 6--9, 1992.
\newblock \doi{10.2514/6.1992-439}.

\bibitem[Stuart and Teckentrup(2017)]{Stuart2017}
A.~M. Stuart and A.~L. Teckentrup.
\newblock Posterior consistency for {Gaussian} process approximations of
  {Bayesian} posterior distributions.
\newblock \emph{Mathematics of Computation}, 87\penalty0 (310):\penalty0
  721--753, 2017.
\newblock \doi{10.1090/mcom/3244}.

\bibitem[Sullivan(2015)]{Sullivan2015}
T.~J. Sullivan.
\newblock \emph{Introduction to Uncertainty Quantification}.
\newblock Springer International Publishing, 2015.
\newblock \doi{10.1007/978-3-319-23395-6}.

\bibitem[Tokdar and Kass(2009)]{Tokdar2009}
S.~T. Tokdar and R.~E. Kass.
\newblock Importance sampling: a review.
\newblock \emph{Wiley Interdisciplinary Reviews: Computational Statistics},
  2\penalty0 (1):\penalty0 54--60, 2009.
\newblock \doi{10.1002/wics.56}.

\bibitem[Vrugt et~al.(2009)Vrugt, ter Braak, Diks, Robinson, Hyman, and
  Higdon]{Vrugt2009}
J.~A. Vrugt, C.~J.~F. ter Braak, C.~G.~H. Diks, B.~A. Robinson, J.~M. Hyman,
  and D.~Higdon.
\newblock Accelerating {Markov} chain {Monte} {Carlo} simulation by
  differential evolution with self-adaptive randomized subspace sampling.
\newblock \emph{International Journal of Nonlinear Sciences and Numerical
  Simulation}, 10\penalty0 (3):\penalty0 273--290, 2009.
\newblock \doi{10.1515/ijnsns.2009.10.3.273}.

\bibitem[Watson(1980)]{Watson1980}
G.~A. Watson.
\newblock \emph{Approximation Theory and Numerical Methods}.
\newblock John Wiley \& Sons Ltd, 1980.
\newblock ISBN 0471277061.

\bibitem[Wilcox(1998)]{Wilcox1998}
D.~C. Wilcox.
\newblock \emph{Turbulence modeling for CFD}, volume~2.
\newblock DCW industries La Ca\~nada, CA, 1998.

\bibitem[Witteveen et~al.(2009)Witteveen, Doostan, Pecnik, and
  Iaccarino]{Witteveen2009a}
J.~A.~S. Witteveen, A.~Doostan, R.~Pecnik, and G.~Iaccarino.
\newblock Uncertainty quantification of the transonic flow around the {RAE
  2822} airfoil.
\newblock \emph{Center for Turbulence Research, Annual Briefs, Stanford
  University}, 2009.

\bibitem[Witteveen and Iaccarino(2010)]{Witteveen2010}
J.~A.~S. Witteveen and G.~Iaccarino.
\newblock Simplex elements stochastic collocation for uncertainty propagation
  in robust design optimization.
\newblock In \emph{48\textsuperscript{th} {AIAA} Aerospace Sciences Meeting
  Including the New Horizons Forum and Aerospace Exposition}. American
  Institute of Aeronautics and Astronautics ({AIAA}), 2010.
\newblock \doi{10.2514/6.2010-1313}.

\bibitem[Witteveen and Iaccarino(2012)]{Witteveen2012}
J.~A.~S. Witteveen and G.~Iaccarino.
\newblock Simplex stochastic collocation with random sampling and extrapolation
  for nonhypercube probability spaces.
\newblock \emph{{SIAM} J. Sci. Comput.}, 34\penalty0 (2):\penalty0 A814--A838,
  2012.
\newblock \doi{10.1137/100817504}.

\bibitem[Witteveen and Iaccarino(2013)]{Witteveen2013}
J.~A.~S. Witteveen and G.~Iaccarino.
\newblock Simplex stochastic collocation with {ENO}-type stencil selection for
  robust uncertainty quantification.
\newblock \emph{Journal of Computational Physics}, 239:\penalty0 1--21, 2013.
\newblock \doi{10.1016/j.jcp.2012.12.030}.

\bibitem[Xiu(2010)]{Xiu2010}
D.~Xiu.
\newblock \emph{Numerical Methods for Stochastic Computations}.
\newblock Princeton University Press, 2010.
\newblock ISBN 0691142122.

\bibitem[Yan and Zhang(2017)]{Yan2017}
L.~Yan and Y.-X. Zhang.
\newblock Convergence analysis of surrogate-based methods for {Bayesian}
  inverse problems.
\newblock \emph{Inverse Problems}, 33\penalty0 (12):\penalty0 125001, 2017.
\newblock \doi{10.1088/1361-6420/aa9417}.

\end{thebibliography}

\end{document}